# High order compact fully-discrete scheme for hyperbolic conversation laws

Tong Zhou[a,c], Haitao Dong[b], Shucheng Pan[a,c,*]

[a] School of Aeronautics, Northwestern Polytechnical University, Xi'an, 710072, PR China

[b] School of aeronautic science and engineering, Beihang University, Beijing, 100191, PR China

[c] National key laboratory of aircraft configuration design, Xi'an, 710072, PR China

**Abstract:** Based on the solution formula method, a series of one-step fully-discrete schemes, such as FWENO/Full-WENO has been proposed. Storing the by-products conservative variables at the half points (grid center) and using them as interpolation information to construct high-order schemes, we obtain a new class of one-step compact fully-discrete schemes. The new scheme can be associate with various non-oscillatory strategies. This paper takes state-of-the-art WENO-JS method as an example and proposes a family of compact fully-discrete WENO scheme. Detailed analysis is conducted on accuracy, errors, computational cost, efficiency and their connection with Hermite interpolation. Meanwhile, we design a new entropy flux linearization strategy for Euler equations to enhance its robustness, and also develop a multi-dimensional method for this compact fully-discrete framework. Due to the new scheme is one-step and utilizes stored by-products information for interpolation, it has a significant advantage in efficiency. For one-dimensional Euler equations, compared to the original FWENO, the computational cost only increases by 20-40%, while is approximately one-third of WENO+RK3. For two-dimensional case, a new special dimension-by-dimension strategy is applied. Although there is an additional computing cost, numerical experiments show that the new scheme only needs about 1/10 to 1/13 cost of that for WENO+RK3 when obtaining similar or even better resolution results, indicating that the new scheme is more efficient than semi-discrete schemes based on RK methods.

## 1 Introduction

As the governing equations for fluid mechanics, Euler equations are hyperbolic conservation laws, which may generate shock discontinuities even if the initial values are sufficiently smooth. This presents a significant challenge in capturing flow structures, making the development of advanced numerical schemes a long-term research topic in computational fluid dynamics. In 1983, Harten [1] introduced the concept of total variation diminishing (TVD) and proposed the second-order TVD scheme. However, the stringent nature of the TVD condition made it challenging to develop to high-order. Then Harten et al. [2–5] proposed the Essentially Non-Oscillatory (ENO) scheme, which uncovered the prelude of high-order schemes. ENO schemes always select the smoothest stencil, which results in the wastage of significant stencil information. Thus, Liu et al. [6] first proposed the Weighted Essentially Non-Oscillatory (WENO) scheme, which weight $r$ stencils form ($r$)-th to ($2r$-1)-th order by smoothness indicators. It makes WENO achieve ($2r$-1)-th order at smooth regions. Soon later, through the efforts of Jiang and Shu [7], the WENO scheme saw significant advancements. Now, a large portion of high-order schemes have been developed based on the ENO/WENO concepts, with improvements either in smoothness indicators or stencil selection strategy, such as WENO-MP [8], WENO-M [9], WENO-Z [10–13], WENO-AO [14], WENO-CU [15], TENO [16–20], WENO-MR [21–24], and so on. Currently, widely used high-accuracy schemes include the Discontinuous Galerkin (DG) [25–31] method, Weighted Compact Nonlinear Scheme (WCNS) [32–34], Hermite WENO (HWENO)[35–38].

[*] Corresponding author.
E-mail address: shucheng.pan@nwpu.edu.cn (S.Pan).

The aforementioned high-order schemes are all discretized in a semi-discrete form, so they can also be referred to as semi-discrete schemes. These semi-discrete schemes primarily only discretize the spatial derivatives, while the temporal derivatives are generally discretized by Runge-Kutta (RK) method, transforming the partial differential equations into ordinary differential equations. Although the explicit RK method has been fully developed at present, there still exist two barriers difficult to solve. (1) Accuracy barrier: RK method with higher than fourth orders cannot guarantee the TVD property. Therefore, in most situations, we can only apply TVD-RK3 [7] for time evolution; (2) Efficiency barrier: For an explicit $s$-stage $p$-th order RK method, there is $s \geq p$, and $s > p$ when $p \geq 5$. The number of stage $s$ in RK method corresponds to the times of global spatial discretization required for each time step. For the widely used 3-stage TVD-RK3, the computational cost is three times that of one-step schemes, which significantly limits computational efficiency. Also, the memory requirement is three times that of one-step schemes. At the same time, accuracy and efficiency barriers are interrelated: relatively low-order RK methods and smaller CFL number can be used to balance computational cost and accuracy.

Of course, in practical computations, spatial accuracy often has a much greater impact on the results than temporal accuracy, which results in TVD-RK3 can also yield satisfactory results. Additionally, designing high-order fully-discrete schemes with coupled space-time discretization is very challenging. So, most researchers prefer to design semi-discrete schemes based on RK methods. However, for long-term evolution problems such as compressible turbulence, both spatial and temporal accuracy are crucial. Consequently, some researchers are now focusing on fully discrete schemes. However, the development of fully-discrete schemes lags behind that of semi-discrete schemes, due to their complexity and difficulty.

Till now, most fully-discrete schemes are based on the Lax-Wendroff (LW) method [39–41], which transforms the temporal derivatives into spatial derivatives through Taylor expansion, enabling the simultaneous discretization of both spatial and temporal derivatives. With the Lax-Wendroff framework, successful schemes include the Arbitrary High Order Derivative Riemann problem (ADER) [42–45] approach and the Lax-Wendroff type WENO/HWENO [40,46,47] schemes are proposed. ADER scheme combines the LW method with the Generalized Riemann Problem (GRP) solver to achieve a one-step fully discrete scheme. For linear constant coefficient equation, ADER scheme can be 2-3 times faster than the ENO/WENO schemes. However, for the one-dimensional Euler equations, ADER3 is only 50% faster than WENO5 [42]. For two-dimensional Euler equations, ADER3 schemes are faster than the WENO scheme roughly by 70% [43], ADER4 in two space dimensions is more expensive than ADER3 by a factor of three. As to the memory requirement, ADER schemes of any order effectively need only two global arrays to store the vector of the conservative variables and the total sum of fluxes. The WENO schemes with the 3-stage TVD-RK3 need at least three such arrays. LW-type WENO/HWENO schemes experience a significant increase in computational cost as accuracy and equation complexity rise, due to the need of computing numerous high-order derivatives. This severely impacts efficiency, which is why in current literatures, one-step fully-discrete schemes based on the LW method rarely exceed fourth order. Generally, multi-stage methods are often required to ensure computational efficiency to a certain extent [48–50].

In summary, there is still significant room for improvement in the efficiency of fully-discrete schemes. Achieving efficient arbitrary high-order one-step fully-discrete schemes, while maintaining high efficiency for systems and multi-dimensional cases, is a common goal for researchers. In recent years, based on the solution formula method proposed by DONG [51] in 2002, DONG&ZHOU developed a new high-order one-step fully-discrete framework, and constructed entropy condition (EC) schemes [52–54] and fully-discrete WENO (FWENO/Full-WENO) schemes [55,56]. The idea behind the solution formula method is to construct and discretize the solution of equation to obtain a conservative scheme. However, the solution of hyperbolic conservation law is difficult to obtain directly. Therefore, the solution formula method integrates the conservation law equations once in space, resulting in the corresponding Hamilton-Jacobi (HJ) equations. The HJ equation has an exact solution for linear fluxes. By linearizing the flux, one can obtain its (quasi) exact solutions, thus providing the expression for numerical flux. The numerical flux of the HJ equation is formally the same as that of the conservation law, thereby enabling the derivation of the corresponding conservative scheme. Fully-discrete schemes within this framework can achieve arbitrary space-time accuracy in one step and offer high computational efficiency. Under same computational conditions, the computational speed of such fully-discrete schemes can roughly reach $s$ times that of semi-discrete schemes using an $s$-stage RK method. For example, FWENO can be approximately three times faster than 3-stage TVD-RK3-WENO.

In this paper, we find that in the fully-discrete framework based on the solution formula method, not only $t^{n+1}$ moment conservative variables at the grid points can be obtained, but also that at the half points. If we store the conservative variables at half points and use them for interpolation, we can construct compact schemes. This allows us to achieve higher accuracy order within same grid dependency range. Moreover, because in the fully-discrete framework based on the solution formula method, the conservative variables at half points are already computed, it does not significantly increase computational cost, at least in the one-dimensional case. Based on this approach and combined with any non-oscillatory interpolation method, we can construct specific schemes. In this paper, we apply state-of-the-art WENO reconstruction as an example to construct the Compact Fully-Discrete WENO schemes (CFWENO). The paper provides the specific forms of CFWENO3, CFWENO5, and CFWENO7. For multi-dimensional cases, Strang splitting [57] in a dimension-by-dimension manner leads to significant errors because it evolves the temporal process in only one direction at the half points. Therefore, we design a new approach to enhance accuracy at half points in these scenarios. Furthermore, this paper also conducted a comprehensive analysis of the compact fully discrete framework based on the solution formula method, covering aspects such as accuracy, errors, computational cost, computational efficiency, and the relation between Hermite interpolation. Extensive numerical experiments have been carried out, which indicate CFWENO has higher computational efficiency than FWENO and WENO+RK3. For the two-dimensional case, when achieving similar or better results, the computational cost for CFWENO is only about one-third of FWENO and about one-tenth of WENO+RK3.

The rest of this article is as follows: The second section briefly reviews semi-discrete schemes based on RK methods and fully-discrete schemes based on the solution formula method; The third section presents the construction process of CFWENO, and extends it to the Euler equations and multi-dimensional cases. It includes a comprehensive analysis of accuracy, errors, the relationship with Hermite interpolation, computational cost, and computational efficiency; The fourth section gives numerical experiments, including tests on accuracy, single equations, and one-dimensional and two-dimensional Euler equations. We also compare the time cost to verify its efficiency; The fifth section provides a summary.

## 2 Review

In this section, we briefly review semi-discrete schemes that use RK methods for time evolution, as well as fully discrete frameworks based on the solution formula method, which achieve high order in both space and time without requiring RK methods or LW approach.

### 2.1 Brief review of semi-discrete schemes via RK method

Considering the initial value problem for a scalar conservation law

$$\begin{cases} u_t + f(u)_x = 0, & (x,t) \in R\times(0,T] \\ u(x,0) = u_0(x), & x \in R \end{cases}. \tag{1}$$

For convenience, we divide the uniform grid into several cells $I_i = [x_{i-1/2}, x_{i+1/2}]$, denote grid step length as $x_{i+1/2} - x_{i-1/2} = h$, $x_i$ is grid node. Integrating Eq.(1) over the cell $I_i$, we obtain the following integral form.

$$\frac{d\overline{u}_i(t)}{dt} = -\frac{1}{h}\left( f(u(x_{i+1/2},t)) - f(u(x_{i-1/2},t)) \right). \tag{2}$$

Then we can derive the semi-discrete scheme

$$\frac{d\overline{u}_i(t)}{dt} = L(u)_i \approx -\frac{1}{h}\left( \tilde{f}_{i+1/2} - \tilde{f}_{i-1/2} \right), \tag{3}$$

The numerical flux is defined as follows.

$$\tilde{f}_{i+1/2} = \tilde{f}^{Riemann}\left( u^-_{i+1/2}, u^+_{i+1/2} \right), \tag{4}$$

$\tilde{f}^{Riemann}\left( u^-_{i+1/2}, u^+_{i+1/2} \right)$ is a Riemann solver, such as Lax-Friedrichs (LF) flux, Roe flux, HLL flux, HLLC flux, etc. Where $u^{\pm}_{i+1/2}$ is the approximate value of $u(x)$ at $x = x_{i+1/2}$. WENO schemes apply following nonlinear weight strategy, weighting $r$ stencils from ($r$)-th order to (2$r$-1)-th order.

$$u^{\pm}_{i+1/2} = \sum_{k=0}^{r-1} \omega_k p_k(x_{i+1/2}), \quad \omega_k = \frac{\alpha_k}{\sum_{k=0}^{r-1}\alpha_k}, \quad \alpha_k = \frac{\gamma_k}{(\beta_k+\varepsilon)^2}, \quad \beta_k = \sum_{l=1}^{r-1}\int_{x_{i-1/2}}^{x_{i+1/2}} h^{2l-1}\left( p_k^{(l)}(x) \right)^2 dx, \tag{5}$$

Where $p(x)$ is a polynomial approximation. Then, Eq.(1) can be transformed into an ordinary differential equation and discretized in time by RK method. The most commonly used method is the 3-stage TVD-RK3, given as follows.

$$\begin{cases} u_t = L(u), \\ L_1 \equiv L(u^n), \\ L_2 \equiv L(u^n + \tau L_1), \\ L_3 \equiv L(u^n + \tau(\frac{1}{4}L_1 + \frac{1}{4}L_2)), \\ u^{n+1} = u^n + \tau(\frac{1}{6}L_1 + \frac{1}{6}L_2 + \frac{4}{6}L_3), \end{cases} \tag{6}$$

Where $\tau = \Delta t$ is the time step.

### 2.2 Brief review of fully-discrete framework via solution formula method

The solution formula method constructs a conservative scheme for hyperbolic conservation laws by constructing and discretizing the solutions of the Hamilton-Jacobi (HJ) equations. The construction process includes following steps: (1) deriving the HJ equation from the hyperbolic conservation law; (2) linearizing the flux of the HJ equation to obtain its (quasi-)exact solution; (3) discretizing the HJ equation to obtain the numerical flux expression of hyperbolic conservation law.

**2.2.1 Fully-discrete framework via solution formula method**

Consider scalar conservation law Eq.(1), integrating once along the spatial direction $x$, and let $\upsilon \equiv \int u dx$, then we can obtain HJ equation

$$\begin{cases} \upsilon_t + f(\upsilon_x) = 0, & (x,t) \in R\times(0,T] \\ \upsilon(x,0) = \upsilon_0(x), & x \in R \end{cases}. \tag{7}$$

Rewrite the flux $f(u)$ into following linear form

$$f(u) = au - f^*, \tag{8}$$

where $a=f_u(u^{n+1})$ is eigenvalue and $f^*=au^{n+1}-f(u^{n+1})$ is local constant, $u^{n+1}$ is the value of $u$ at $t^{n+1}$, which will be given later. Then the solution of HJ equation Eq.(7) can be given as follows

$$\upsilon(x,t) = \upsilon_0(x-at) + tf^*. \tag{9}$$

Eq.(9) can be written in the following discrete form

$$\upsilon_i^{n+1} = \upsilon^n(x_i - a\tau) + \tau f^*. \tag{10}$$

By discretizing the HJ equation and its solution, the expression of numerical flux can be obtained

$$\tilde{f}_{i+1/2} = \frac{\upsilon_{i+1/2}^n - \upsilon_{i+1/2}^{n+1}}{\tau} = \frac{\upsilon_{i+1/2}^n - \upsilon^n(x_{i+1/2} - a\tau)}{\tau} - f^*, \tag{11}$$

where $\tau=\Delta t$ is temporal step, and $h=\Delta x$ is spatial step length. Eq.(11) can be written into following form

$$\begin{cases} \tilde{f}_{i+1/2} = a\bar{u}_{i+1/2} - f^*, \\ \bar{u}_{i+1/2} \equiv \dfrac{\upsilon_{i+1/2}^n - \upsilon^n(x_{i+1/2} - a\tau)}{a\tau} = \dfrac{1}{a\tau}\displaystyle\int_{x_{i+1/2}-a\tau}^{x_{i+1/2}} u^n(x)dx. \end{cases} \tag{12}$$

where $\bar{u}_{i+1/2}$ can be obtained from the interpolation initial value function $v^n(x)$ of HJ equation with Eq.(12), thus can be named initial value reconstruction. Local constants $a$ and $f^*$ are required for flux linearization, so the process of solving $a$ and $f^*$ can be called flux (linearization) reconstruction. $a$ and $f^*$ can be expressed as the form of $u_{i+1/2}^{n+1}$, and $u^{n+1}(x)$ can be express as $\upsilon_x^{n+1}(x)$. Then we have following derivation

$$\begin{cases} a = f_u(u_{i+1/2}^{n+1}), \\ f^* = au_{i+1/2}^{n+1} - f(u_{i+1/2}^{n+1}), \\ u^{n+1}(x) = \upsilon_x^{n+1}(x) = \upsilon_x^n(x - a\tau) = u^n(x - a\tau). \end{cases} \tag{13}$$

Where eigenvalue $a$ in Eq.(13) can be get from $u_{i+1/2}^{n+1}$, while there is $a$ exists in the expression of $u_{i+1/2}^{n+1}$. So, it is an implicit expression, that need to be solved iteratively. However, the iteration does not significantly affect the computational results. In general, multiple iterations are not needed, except for accuracy test, see section 3.2.1 and section 4.1.2. Eq.(12) and Eq.(13) show that we can obtain $\bar{u}_{i+1/2}$ and $u_{i+1/2}^{n+1}$ for final numerical flux as long as we get the expression of $u^n(x)$. The numerical flux obtained from this framework inherently has consistent high-order temporal and spatial accuracy, so the resulting scheme is fully-discrete in one-step. Thus, the numerical flux and the corresponding conservative scheme of the conservation law Eq.(1) are obtained.

$$u_i^{n+1} = u_i^n - \frac{\tau}{h}\left(\tilde{f}_{i+1/2} - \tilde{f}_{i-1/2}\right). \tag{14}$$

Furthermore, within the framework of solution formula method, for linear flux, $a$ and $f^*$ are exact, so the accuracy only depends on initial value reconstruction. Currently, total accuracy of scheme = accuracy of initial value reconstruction. For nonlinear flux, the accuracy depends on both of initial value reconstruction and flux reconstruction. From the error analysis in [52–56], total accuracy of scheme = min (accuracy of initial value reconstruction, 2×accuracy of flux reconstruction).

# 3 Construction of high order compact fully-discrete scheme

The solution formula method not only provides the node values $u_i^{n+1}$ at $t^{n+1}$ in Eq.(14), but also the half points values $u_{i+1/2}^{n+1}$ at $t^{n+1}$ in Eq.(13). If we store the half points values $u_{i+1/2}^{n+1}$ for the interpolation of next moment, we can construct a higher-order compact scheme within same dependency domain. Additionally, since $u_{i+1/2}^{n+1}$ needs to be computed in original fully-discrete schemes, we do not need to add extra equations. Thus, the increase in computational cost is not significant and depends mainly on the complexity of interpolation. This is true for one-dimensional cases, while multi-dimensional cases will be discussed further in section 3.3. In terms of the non-oscillation method, most of existing strategies can be used, such as linear upwind scheme, ENO, WENO, TENO, WENO-Z, WENO-CU, WENO-MR. In this paper, we apply the widely-used WENO-JS reconstruction, weight $r$ stencils from ($r$)-th to (2$r$-1)-th, resulting in the compact fully-discrete WENO scheme (CFWENO). This section uses third, fifth, and seventh-order schemes ($r$=2,3,4) as examples to introduce the construction process of CFWENO based on the solution formula method.

## 3.1 Compact Fully-discrete WENO under solution formula method framework

The conclusion from section 2.2 indicates that the construction of fully-discrete scheme based on solution formula framework requires two tasks: (1) initial value reconstruction, which involves interpolating the initial value function $v(x)$ of HJ equation to obtain the $\bar{u}_{i+1/2}$ in Eq.(12); (2) flux reconstruction, which involves introducing the local constants $a$ and $f^*$ through the flux linearization when constructing the solution of HJ equation.

Section 2.2 shows that we can get all the variables from Eq.(12) and Eq.(13) once we obtain $u^n(x)$. Different from the semi-discrete schemes need the values at both sides of interface for Riemann solver, the fully-discrete framework based on solution formula method uses upwind interpolation according to the sign of eigenvalues, as follows

$$\bar{u}_{i+1/2} = \begin{cases} \bar{u}_{i+1/2}^{-} & a \ge 0 \\ \bar{u}_{i+1/2}^{+} & a \le 0 \end{cases}, \quad u_{i+\frac{1}{2}}^{n+1} = \begin{cases} u_{i+1/2}^{-,n+1} & a \ge 0 \\ u_{i+1/2}^{+,n+1} & a \le 0 \end{cases}. \tag{15}$$

Following text will introduce $\bar{u}_{i+1/2}^{-}$ and $u_{i+1/2}^{-,n+1}$ (the superscript "-" will be omitted), $\bar{u}_{i+1/2}^{+}$ and $u_{i+1/2}^{+,n+1}$ can be obtained by symmetry properties. In this paper, we apply the WENO-JS reconstruction.

As mentioned above, using a corresponding polynomial $p(x)$ to approximate $u(x)$ yields a corresponding compact stencil. For the approximate of CFWENO3, we choose spatial stencils $S$={$I_{i\text{-}1/2}$, $I_i$, $I_{i\text{-}1/2}$}, we can construct two linear polynomials $p_k(x)$($k$=0,1) (for two 2-th order sub-stencil) and a quadratic polynomial $q(x)$ (for the 3-th order big stencil).

$$
\begin{cases}
\frac{1}{h}\int_{x_{j-1/2}}^{x_{j+1/2}} p_0(x)dx = u_j\,, j = i, & p_0(x_j) = u_j, j = i-\frac{1}{2}, \\
\frac{1}{h}\int_{x_{j-1/2}}^{x_{j+1/2}} p_1(x)dx = u_j\,, j = i, & p_1(x_j) = u_j, j = i+\frac{1}{2}, \\
\frac{1}{h}\int_{x_{j-1/2}}^{x_{j+1/2}} q(x)dx = u_j\,, j = i, & q(x_j) = u_j, j = i-\frac{1}{2}, i+\frac{1}{2}.
\end{cases}
\tag{16}
$$

For the approximate of CFWENO5, we choose spatial stencils $S=\{I_{i-1}, I_{i-1/2}, I_i, I_{i-1/2}, I_{i+1}\}$, we can construct three quadratic polynomials $p_k(x)(k=0,1,2)$ (for three 3-th order sub-stencil) and a quartic polynomial $q(x)$ (for the 5-th order big stencil).

$$
\begin{cases}
\frac{1}{h}\int_{x_{j-1/2}}^{x_{j+1/2}} p_0(x)dx = u_j\,, j = i-1, i, & p_0(x_j) = u_j, j = i-\frac{1}{2}, \\
\frac{1}{h}\int_{x_{j-1/2}}^{x_{j+1/2}} p_1(x)dx = u_j\,, j = i, & p_1(x_j) = u_j, j = i-\frac{1}{2}, i+\frac{1}{2}, \\
\frac{1}{h}\int_{x_{j-1/2}}^{x_{j+1/2}} p_2(x)dx = u_j\,, j = i, i+1, & p_2(x_j) = u_j, j = i+\frac{1}{2}, \\
\frac{1}{h}\int_{x_{j-1/2}}^{x_{j+1/2}} q(x)dx = u_j\,, j = i-1, i, i+1, & q(x_j) = u_j, j = i-\frac{1}{2}, i+\frac{1}{2}.
\end{cases}
\tag{17}
$$

For the approximate of CFWENO7, we choose spatial stencils $S=\{I_{i-3/2}, I_{i-1}, I_{i-1/2}, I_i, I_{i-1/2}, I_{i+1}, I_{i+3/2}\}$, we can construct four cubic polynomials $p_k(x)(k=0,1,2,3)$ (for 4-th order sub-stencil) and a sixth-degree polynomial $q(x)$ (for the 7-th order big stencil).

$$
\begin{cases}
\frac{1}{h}\int_{x_{j-1/2}}^{x_{j+1/2}} p_0(x)dx = u_j\,, j = i-1, i, & p_0(x_j) = u_j, j = i-\frac{3}{2}, i-\frac{1}{2}, \\
\frac{1}{h}\int_{x_{j-1/2}}^{x_{j+1/2}} p_1(x)dx = u_j\,, j = i-1, i, & p_1(x_j) = u_j, j = i-\frac{1}{2}, i+\frac{1}{2}, \\
\frac{1}{h}\int_{x_{j-1/2}}^{x_{j+1/2}} p_2(x)dx = u_j\,, j = i, i+1, & p_2(x_j) = u_j, j = i-\frac{1}{2}, i+\frac{1}{2}, \\
\frac{1}{h}\int_{x_{j-1/2}}^{x_{j+1/2}} p_3(x)dx = u_j\,, j = i, i+1, & p_3(x_j) = u_j, j = i+\frac{1}{2}, i+\frac{3}{2}, \\
\frac{1}{h}\int_{x_{j-1/2}}^{x_{j+1/2}} q(x)dx = u_j\,, j = i-1, i, i+1, & q(x_j) = u_j, j = i-\frac{3}{2}, i-\frac{1}{2}, i+\frac{1}{2}, i+\frac{3}{2}.
\end{cases}
\tag{18}
$$

Accordingly, the strategy for CFWENO is shown in Fig. 1, where circles represent the values $u_j$ ($j=\ldots, i-1, i, i+1, \ldots$) at the grid points obtained from Eq.(14), and diamonds represent the cell-centered values $u_j$ ($j=\ldots i-1/2, i+1/2\ldots$) obtained from Eq.(13)

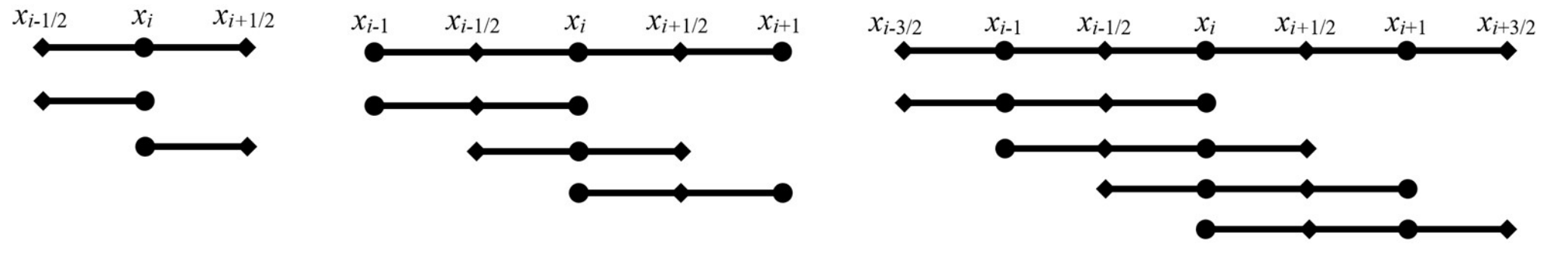

(A) CFWENO3 (B) CFWENO5 (C) CFWENO7

**Fig. 1 Sketch of reconstruction schemes: (A) CFWENO3, (B) CFWENO5, (C) CFWENO7.**

In Eq.(16)(17)(18), we approximate $u_j$ ($j=\ldots, i-1, i, i+1,\ldots$) by the integral average of $p(x)$ in $[x_{j-1/2}, x_{j+1/2}]$, while directly approximate $u_j$($j=\ldots, i-1/2, i+1/2, \ldots$) by $p(x_j)$. This is because $\overline{u}_{i+1/2}$ is an integral average value, and $u_{i+1/2}^{n+1}$ is obtained directly by interpolation.

### 3.1.2 Initial value reconstruction for compact fully-discrete WENO

With the polynomial $p(x)$ obtained in section 3.1.2, we can easily derive the $\bar{u}_{i+1/2}$ from Eq.(12) for initial value reconstruction. The specific form is detailed in Appendix A

$$\bar{u}_{i+1/2}=\frac{1}{vh}\int_{x_{i+1/2}-vh}^{x_{i+1/2}} u^n(x)dx \approx \frac{1}{vh}\int_{x_{i+1/2}-vh}^{x_{i+1/2}} q(x)dx, \tag{19}$$

Similarly, we can now derive the interpolated values $u_{i+1/2}^{n+1}$ needed for the initial value reconstruction at the next time step $t^{n+1}$ from Eq.(13). These values will also be used for the flux reconstruction at the current time $t^n$. The specific form is detailed in Appendix A.

$$u_{i+1/2}^{n+1}= u^n(x_{i+1/2}-vh) \approx q(x_{i+1/2}-vh). \tag{20}$$

Based on the WENO reconstruction concept, the $(2r-1)$-th order polynomial $q(x)$ can be weighted by $(r)$-th order polynomial $p_k(x)$ ($k$=0, 1, 2,…,$r$-1)

$$q(x)=\sum_{k=0}^{r-1}\gamma_k^r p_k^r(x), \tag{21}$$

where $\gamma_k^r$ ($k$=0,1,…,$r$-1) is the optimal linear weight. To achieve order reduction at discontinuities and ensure non-oscillatory properties, this paper employs WENO-JS reconstruction.

$$\bar{u}_{i+1/2}^{2r-1} \approx \sum_{k=0}^{r-1}\bar{\omega}_k^r \frac{1}{vh}\int_{x_{i+1/2}-vh}^{x_{i+1/2}} p_k^r(x)dx, \quad u_{i+1/2}^{n+1,2r-1} \approx \sum_{k=0}^{r-1}\omega_k^r p_k^r(x_{i+1/2}-vh), \tag{22}$$

where, $\bar{\omega}_k^r$ and $\omega_k^r$ are nonlinear weight for $\bar{u}_{i+1/2}$ and $u_{i+1/2}^{n+1}$ respectively.

$$\bar{\omega}_k^r=\frac{\bar{\alpha}_k^r}{\sum_{k=0}^{r-1}\bar{\alpha}_k^r}, \quad \bar{\alpha}_k=\frac{\bar{\gamma}_k^r}{(\beta_k^r+\varepsilon)^2}, \quad \omega_k^r=\frac{\alpha_k^r}{\sum_{k=0}^{r-1}\alpha_k^r}, \quad \alpha_k^r=\frac{\gamma_k^r}{(\beta_k^r+\varepsilon)^2}. \tag{23}$$

The optimal linear weight $\bar{\gamma}_k^r$ ($k$=0,1,…,$r$-1) for $\bar{u}_{i+1/2}$ (the superscript $(2r-1)$ will be omitted) is given in Table 1

**Table 1 Optimal linear weights $\bar{\gamma}_k^r$ for $\bar{u}_{i+1/2}$**

| $\bar{\gamma}_k^r$ | $k=0$ | $k=1$ | $k=2$ | $k=3$ |
|---|---|---|---|---|
| $r=2$ | $v$ | $1-v$ | —— | —— |
| $r=3$ | $\frac{1}{6}v(1+v)$ | $\frac{1}{6}(1+v)(2-v)$ | $\frac{1}{6}(1-v)(2-v)$ | —— |
| $r=4$ | $\frac{1}{18}v(1+v)^2$ | $\frac{7}{54}(1+v)^2(2-v)$ | $\frac{7}{54}(1+v)(2-v)^2$ | $\frac{1}{18}(1-v)(2-v)^2$ |

The optimal linear weight $\gamma_k^r$ ($k$=0,1,…,$r$-1) for $u_{i+1/2}^{n+1}$ (the superscript $(2r-1)$ will be omitted) is given in Table 2

**Table 2 Optimal linear weights $\gamma_k^r$ for $u_{i+1/2}^{n+1}$**

| $\gamma_k^r$ | $k=0$ | $k=1$ | $k=2$ | $k=3$ |
|---|---|---|---|---|
| $r=2$ | $\frac{3v^2-2v}{2v-1}$ | $\frac{-3v^2+4v-1}{2v-1}$ | —— | —— |
| $r=3$ | $\frac{1}{6}\frac{v(5v^2+v-2)}{3v-1}$ | $-\frac{1}{6}\frac{30v^4-60v^3-v^2+31v-8}{(3v-1)(3v-2)}$ | $\frac{1}{6}\frac{(v-1)(5v^2-11v+4)}{(3v-2)}$ | —— |
| $r=4$ | $\gamma_0^4$ | $\gamma_1^4$ | $\gamma_2^4$ | $\gamma_3^4$ |

where

$$
\begin{aligned}
\gamma_0^4 &= \frac{1}{36}\frac{v(1+v)(7v^3-12v^2-3v+4)}{2v^2-4v+1}, \\
\gamma_1^4 &= \frac{1}{108}\frac{(1+v)(2-v)(98v^4-177v^3+4v^2+61v-14)}{(2v^2-4v+1)(2v-1)}, \\
\gamma_2^4 &= \frac{1}{108}\frac{(1+v)(v-2)(98v^4-215v^3+61v^2+70v-28)}{(2v-1)(2v^2-1)}, \\
\gamma_3^4 &= \frac{1}{36}\frac{(1-v)(v-2)(7v^3-9v^2-6v+4)}{2v^2-1}.
\end{aligned}
\tag{24}
$$

The aforementioned weight for $u_{i+1/2}^{n+1}$ may potentially become negative or result in a zero denominator. For negative weight, we employ the method described in reference [58], which requires little additional computational cost. For the zero denominator problem, we cut off $\gamma_k^r(v)$ by $\gamma_k^r(\delta_k^r-\varepsilon_d)$ or $\gamma_k^r(\delta_k^r+\varepsilon_d)$, according to whether $v$ is at left or right side of the possible discontinuities between $v\in(\delta_k^r-\varepsilon_d,\delta_k^r+\varepsilon_d)$. Where $\delta_k^r$ is the position of discontinuities, see in Table 3. $\varepsilon_d$ is a small positive number satisfy $0<\delta_{k-1}^r<\delta_k^r-\varepsilon_d$ and $\delta_k^r+\varepsilon_d<\delta_{k+1}^r<1$, in this paper we choose $\varepsilon_d=0.05$.

**Table 3 Position of discontinuities $\delta_k^r$ for $\gamma_k^r$ of $u_{i+1/2}^{n+1}$**

| $\delta_k^r$ | $k=1$ | $k=2$ | $k=3$ |
|---|---|---|---|
| $r=2$ | $1/2$ | —— | —— |
| $r=3$ | $1/3$ | $2/3$ | —— |
| $r=4$ | $1-\sqrt{2}/2$ | $1/2$ | $\sqrt{2}/2$ |

$\bar{\beta}_k^r$ and $\beta_k^r$ in Eq.(23) are smooth indicators for WENO reconstruction. Since slight differences in smoothness indicators do not significantly impact the results, for simplicity, we set $\bar{\beta}_k^r=\beta_k^r$ and use the same simplified calculation method, which takes the following form

$$
\beta_k^r=\sum_{l=1}^{r-1}\left(\int_{x_{i-1/2}}^{x_{i+1/2}} h^{2l-1}\left(\left(\frac{1}{vh}\int_{x_{i+1/2}-vh}^{x_{i+1/2}} p_k^r(x)dx\right)^{(l)}\right)^2 dx\right) \tag{25}
$$

Eq.(25) is quite complex, so we simplify it by setting $v$=0. This makes it similar to the semi-discrete WENO scheme and results in the following explicit expression.

$$
\begin{aligned}
&\beta_0^2 = 4(u_j - u_{j-\frac{1}{2}})^2\\
&\beta_1^2 = 4(u_{j+\frac{1}{2}} - u_j)^2\\
&\beta_0^3 = \frac{1}{4}\left(u_{j-1} - 6u_{j-\frac{1}{2}} + 5u_j\right)^2 + \frac{39}{4}\left(u_{j-1} - 2u_{j-\frac{1}{2}} + u_j\right)^2\\
&\beta_1^3 = \left(u_{j-\frac{1}{2}} - u_{j+\frac{1}{2}}\right)^2 + 39\left(u_{j-\frac{1}{2}} - 2u_j + u_{j+\frac{1}{2}}\right)^2\\
&\beta_2^3 = \frac{1}{4}\left(5u_j - 6u_{j+\frac{1}{2}} + u_{j+1}\right)^2 + \frac{39}{4}\left(u_j - 2u_{j+\frac{1}{2}} + u_{j+1}\right)^2\\
&\beta_0^4 = \left(-u_{j-\frac{3}{2}} + 3u_{j-1} - 5u_{j-\frac{1}{2}} + 3u_j\right)^2 + 39\left(-u_{j-\frac{3}{2}} + 3u_{j-1} - 3u_{j-\frac{1}{2}} + u_j\right)^2 + \frac{781}{20}\left(-2u_{j-\frac{3}{2}} + 5u_{j-1} - 4u_{j-\frac{1}{2}} + u_j\right)^2\\
&\beta_1^4 = \left(u_{j-\frac{1}{2}} - u_{j+\frac{1}{2}}\right)^2 + 39\left(u_{j-\frac{1}{2}} - 2u_j + u_{j+\frac{1}{2}}\right)^2 + \frac{781}{20}\left(-u_{j-1} + 4u_{j-\frac{1}{2}} - 5u_j + 2u_{j+\frac{1}{2}}\right)^2\\
&\beta_2^4 = \left(u_{j-\frac{1}{2}} - u_{j+\frac{1}{2}}\right)^2 + 39\left(u_{j-\frac{1}{2}} - 2u_j + u_{j+\frac{1}{2}}\right)^2 + \frac{781}{20}\left(-2u_{j-\frac{1}{2}} + 5u_j - 4u_{j+\frac{1}{2}} + u_{j+1}\right)^2\\
&\beta_3^4 = \left(3u_j - 5u_{j+\frac{1}{2}} + 3u_{j+1} - u_{j+\frac{3}{2}}\right)^2 + 39\left(u_j - 3u_{j+\frac{1}{2}} + 3u_{j+1} - u_{j+\frac{3}{2}}\right)^2 + \frac{781}{20}\left(u_j - 4u_{j+\frac{1}{2}} + 5u_{j+1} - 2u_{j+\frac{3}{2}}\right)^2
\end{aligned}
\tag{26}
$$

To simplify the computation, we also directly use the smoothness indicator from Eq.(26) for weighting $u_{i+1/2}^{n+1}$.

### 3.1.2 Flux reconstruction with entropy condition for compact fully-discrete WENO

Flux reconstruction ultimately requires $a$ and $f^*$ in flux linearization. These parameters significantly impact stability, especially near discontinuities. Therefore, entropy conditions are commonly used for selection: Roe averaging is applied to shock waves, while high-order reconstruction is used for other cases.

$$
a = \begin{cases} \dfrac{f(u_R) - f(u_L)}{u_R - u_L} & v_L > v_R \\ f_u(u_{i+1/2}^{n+1}) & v_L \le v_R \end{cases}, \quad f^* = \begin{cases} a\dfrac{u_L + u_R}{2} - \dfrac{f(u_L) + f(u_R)}{2} & v_L > v_R \\ a u_{i+1/2}^{n+1} - f(u_{i+1/2}^{n+1}) & v_L \le v_R \end{cases}, \quad \begin{cases} v_L \equiv f_u(u_L)\dfrac{\tau}{h} \\ v_R \equiv f_u(u_R)\dfrac{\tau}{h} \end{cases}
\tag{27}
$$

Eq.(27) represents the entropy condition flux linearization strategy for a scalar conversation law. When $u_{i+1/2}^{n+1} = (u_i + u_{i+1})/2$ it is defined as baseline entropy condition reconstruction. Since this fully discrete framework is designed from scalar equation, in case of systems, the characteristic form is generally used. Flux linearization is performed first, followed by projecting the system onto the characteristic directions by eigenvectors. Thus, each equation behaves similarly to the scalar equation.

$$
\begin{cases} \mathbf{u}_t + \mathbf{f}(\mathbf{u})_x = 0 \\ \mathbf{v}_t + \mathbf{f}(\mathbf{v}_x) = 0 \end{cases} \xrightarrow{\text{linearization}} \begin{cases} \mathbf{u}_t + (\mathbf{A}\mathbf{u} - \mathbf{f}^*)_x = 0 \\ \mathbf{v}_t + \mathbf{A}\mathbf{v}_x - \mathbf{f}^* = 0 \end{cases} \xrightarrow{\text{diagonalization}} \begin{cases} \mathbf{L}\mathbf{u}_t + (\Lambda\mathbf{L}\mathbf{u} - \mathbf{L}\mathbf{f}^*)_x = 0 \\ \mathbf{L}\mathbf{v}_t + \Lambda\mathbf{L}\mathbf{v}_x - \mathbf{L}\mathbf{f}^* = 0 \end{cases}
\tag{28}
$$

Where **L** represents the left eigenvectors, and **A** is the Jacobian matrix of flux. In Euler equations, the baseline entropy condition linearization typically uses velocity for total determination: Roe averaging is used for shock waves, while arithmetic averaging is used for the others, as detailed below

$$
\begin{cases}
u=\begin{cases}\dfrac{\sqrt{\rho_L}u_L+\sqrt{\rho_R}u_R}{\sqrt{\rho_L}+\sqrt{\rho_R}} & u_L>u_R\\ u(\dfrac{\mathbf{u}_L+\mathbf{u}_R}{2}) & u_L\le u_R\end{cases},\ H=\begin{cases}\dfrac{\sqrt{\rho_L}H(\mathbf{u}_L)+\sqrt{\rho_R}H(\mathbf{u}_R)}{\sqrt{\rho_L}+\sqrt{\rho_R}} & u_L>u_R\\ H(\dfrac{\mathbf{u}_L+\mathbf{u}_R}{2}) & u_L\le u_R\end{cases},\ \begin{cases}c=\sqrt{(\gamma-1)(H-0.5u^2)}\\ \Lambda=\{\lambda^k\}=\Lambda(u,c)\\ \mathbf{L}=\{\mathbf{L}^k\}=\mathbf{L}(u,c)\end{cases}\\
\mathbf{u}_b^*=\dfrac{\mathbf{u}_L+\mathbf{u}_R}{2},\ \mathbf{f}_b=\begin{cases}\dfrac{\mathbf{f}(\mathbf{u}_R)+\mathbf{f}(\mathbf{u}_R)}{2} & u_L>u_R\\ \mathbf{f}(\dfrac{\mathbf{u}_L+\mathbf{u}_R}{2}) & u_L\le u_R\end{cases},\ \boldsymbol{\varphi}^*=\Lambda\mathbf{L}\mathbf{u}_b^*-\mathbf{L}\mathbf{f}_b \Leftrightarrow \varphi^{*k}=\lambda^k\mathbf{L}^k\mathbf{u}_b^{*k}-\mathbf{L}^k\mathbf{f}_b
\end{cases}
\tag{29}
$$

In Eq.(29), $\mathbf{u}_L\equiv\mathbf{u}_i$, $\mathbf{u}_R\equiv\mathbf{u}_{i+1}$, $\boldsymbol{\varphi}^*\equiv\mathbf{L}\mathbf{f}^*$, $\varphi^{*k}\equiv\mathbf{L}^k\mathbf{f}^{*k}$. To judge shock accurately and to prevent instability caused by using high-order flux reconstruction near shocks, we design following flux reconstruction strategy (for the Euler equations, high-order flux reconstruction only updates eigenvalues $\boldsymbol{\Lambda}$ and local constant $\mathbf{f}^*$, while the eigenvectors **LR** are directly obtained from the baseline entropy flux reconstruction)

$$
\begin{array}{l}
\textit{if}\ (\max(|p_L|,|p_R|)\ge s_1\min(|p_L|,|p_R|)\cup p_Lp_R\le 0)\quad \text{(option 1)}\\
\quad \text{use baseline entropy condition flux reconstruction}\\
\textit{else}\\
\quad \text{get guess middle pressure } p_m\\
\quad \textit{if}\ (\max(|p_L|,|p_R|)< s_2\min(|p_L|,|p_R|))\qquad \text{(option 2)}\\
\quad\quad \varphi^{*1}\equiv\mathbf{L}^1\left(\lambda^1(\mathbf{u}^*)\mathbf{u}^*-\mathbf{f}(\mathbf{u}^*)\right),\quad \varphi^{*2}\equiv\mathbf{L}^2\left(\lambda^2(\mathbf{u}^*)\mathbf{u}^*-\mathbf{f}(\mathbf{u}^*)\right),\ \varphi^{*3}\equiv\mathbf{L}^3\left(\lambda^3(\mathbf{u}^*)\mathbf{u}^*-\mathbf{f}(\mathbf{u}^*)\right)\\
\quad \textit{else}\\
\quad\quad \textit{if}\,(p_L<p_m\cap p_m>p_R)\qquad \text{(option 3)}\\
\quad\quad\quad \varphi^{*1}\equiv\mathbf{L}^1(\lambda_b^1\mathbf{u}_b^*-\mathbf{f}_b),\quad \varphi^{*2}\equiv\mathbf{L}^2\left(\lambda_b^2\mathbf{u}_b^*-\mathbf{f}_b\right),\quad \varphi^{*3}\equiv\mathbf{L}^3(\lambda_b^3\mathbf{u}_b^*-\mathbf{f}_b)\\
\quad\quad \textit{elseif}\,(p_L\ge p_m\cap p_m\le p_R)\qquad \text{(option 4)}\\
\quad\quad\quad \varphi^{*1}\equiv\mathbf{L}^1\left(\lambda^1(\mathbf{u}^*)\mathbf{u}^*-\mathbf{f}(\mathbf{u}^*)\right),\quad \varphi^{*2}\equiv\mathbf{L}^2\left(\lambda^2(\mathbf{u}^*)\mathbf{u}^*-\mathbf{f}(\mathbf{u}^*)\right),\quad \varphi^{*3}\equiv\mathbf{L}^3\left(\lambda^3(\mathbf{u}^*)\mathbf{u}^*-\mathbf{f}(\mathbf{u}^*)\right)\\
\quad\quad \textit{elseif}\,(p_L<p_m\cap p_m\le p_R)\qquad \text{(option 5)}\\
\quad\quad\quad \varphi^{*1}\equiv\mathbf{L}^1(\lambda_b^1\mathbf{u}_b^*-\mathbf{f}_b),\quad \varphi^{*2}\equiv\mathbf{L}^2\left(\lambda^2(\mathbf{u}^*)\mathbf{u}^*-\mathbf{f}(\mathbf{u}^*)\right),\quad \varphi^{*3}\equiv\mathbf{L}^3\left(\lambda^3(\mathbf{u}^*)\mathbf{u}^*-\mathbf{f}(\mathbf{u}^*)\right)\\
\quad\quad \textit{elseif}\,(p_L\ge p_m\cap p_m>p_R)\qquad \text{(option 6)}\\
\quad\quad\quad \varphi^{*1}\equiv\mathbf{L}^1\left(\lambda^1(\mathbf{u}^*)\mathbf{u}^*-\mathbf{f}(\mathbf{u}^*)\right),\quad \varphi^{*2}\equiv\mathbf{L}^2\left(\lambda^2(\mathbf{u}^*)\mathbf{u}^*-\mathbf{f}(\mathbf{u}^*)\right),\quad \varphi^{*3}\equiv\mathbf{L}^3(\lambda_b^3\mathbf{u}_b^*-\mathbf{f}_b)\\
\quad\quad \textit{endif}\\
\quad \textit{endif}\\
\textit{endif}
\end{array}
\tag{30}
$$

where $\mathbf{u}^*\equiv\mathbf{u}_{i+1/2}^{n+1}$. $p_m$ is estimated interface pressure for judging whether [$\mathbf{u}_L$, $\mathbf{u}_R$] will evolve into a shock wave, given as follows

$$
p_m=\max\left(\min\left((\max(0,\frac{\frac{u_L-u_R}{2}+\frac{c_L}{\gamma-1}+\frac{c_R}{\gamma-1}}{\frac{c_L}{\gamma-1}(\frac{1}{p_L})^{\frac{\gamma-1}{2\gamma}}+\frac{c_R}{\gamma-1}(\frac{1}{p_R})^{\frac{\gamma-1}{2\gamma}}}))^{\frac{2\gamma}{\gamma-1}},(\frac{4}{\frac{1}{\sqrt{p_L}}+\frac{1}{\sqrt{p_R}}})^2\right),(\max(0,\frac{u_L-u_R}{\frac{1}{\sqrt{\rho_L\frac{\gamma+1}{2}}}+\frac{1}{\sqrt{\rho_R\frac{\gamma+1}{2}}}}))^2\right) \tag{31}
$$

The concept of entropy condition flux reconstruction in Eq.(30) is as follows: (1) Option 1 is strong shock and positivity constraint, all use the baseline entropy condition flux reconstruction. Where $s_1\ge1$, the smaller $s_1$ is, the more stable the method. The larger $s_1$ is, the higher the accuracy. In this paper we take $s_1=2$; (2) Option 2 is designed to prevent the use of lower-order reconstruction due to small numerical errors. Where $s_1\ge1$ and close to 1, in this paper we take $s_2$=1.05; (3) Option 3 represents both the left and right waves are shock waves, we all use baseline reconstruction; (4) Option 4 represents both the left and right waves are rarefaction wave, we all use high-order reconstruction; (5) Option 5 represents a left shock wave and a right rarefaction wave, we apply baseline reconstruction for left wave and high-order reconstruction for others; (6) Option 6 represents a left rarefaction wave and a right shock wave, we apply baseline reconstruction for right wave and high-order reconstruction for others. With all above information, we can combine them into the final numerical flux.

$$\tilde{\mathbf{f}}_{j+\frac{1}{2}} = \mathbf{L}^{-1}\{\mathbf{L}^k(\lambda^k \bar{\mathbf{u}}_{j+\frac{1}{2}} - \mathbf{f}^{*k})\}\ , \tag{32}$$

Then, we can obtain the fully-discrete scheme for Euler equations based on solution formula framework.

$$\mathbf{u}_j^{n+1} = \mathbf{u}_j^n - \frac{\tau}{h}(\tilde{\mathbf{f}}_{j+\frac{1}{2}} - \tilde{\mathbf{f}}_{j-\frac{1}{2}})\ . \tag{33}$$

## 3.2 Some analysis for compact fully-discrete WENO

This section mainly analyzes some properties of CFWENO, including accuracy, error and its relationship with Hermite interpolation.

### 3.2.1 Accuracy analysis for compact fully-discrete WENO

For fully-discrete scheme under solution formula method, there are the conclusions from the error analysis in [52–56]

(1) linear accuracy order = accuracy order of initial value reconstruction;

(2) nonlinear accuracy order = 2×accuracy order of flux reconstruction;

(3) total accuracy order = min (accuracy order of initial value reconstruction, 2×accuracy order of flux reconstruction).

Since eigenvalue $a(u_{i+1/2}^{n+1})$ required for initial value reconstruction contains $a$ in $u_{i+1/2}^{n+1} = u^n(x - a\tau)$ . Thus, iteration is required to achieve desired accuracy. For FWENO, each iteration of $a$ improves the accuracy of $u_{i+1/2}^{n+1}$ by one order and improves the accuracy of flux reconstruction by two order (first-order $a$ can meet the second-order accuracy requirement for $u_{i+1/2}^{n+1}$, and so on). However, for CFWENO, $u_{i+1/2}^{n+1}$ is used not only for flux reconstruction but also for initial value reconstruction. For initial value reconstruction, each iteration only increases the accuracy by one order. This phenomenon will be shown in numerical experiment in section 4.1.2. Iterations significantly impact efficiency. Fortunately, actual computations show that, especially for Euler equations, iterations show little impact on results. So, in the computation, we only need to apply the $a$ computed by baseline flux reconstruction to compute the high-order one. All the computational results in this paper used this strategy, except for those involving accuracy test with sufficient iterations.

### 3.2.2 Error analysis for compact fully-discrete WENO

Table 4 and Fig. 2 present the error coefficients $C_e$ for CFWENO, FWENO, and WENO. The error can be given as $R^{2r-1} = C_e u^{(2r-1)} h^{2r-1}$, where $v$=$a\tau/h$≈CFL (for linear flux $v$=CFL). The error coefficient for nonlinear flux is more complex, which needs to account for the error in $v$. This section only discusses the linear case, the nonlinear case is discussed through numerical experiments.

In Fig. 2, it can be observed that the error coefficient for WENO is constant, while that for FWENO decreases with $v$. When $v$=0, the error coefficient of FWENO is the same as that for semi-discrete WENO (since the fully-discrete stencil for FWENO degenerates to the semi-discrete form when $v$=0), when $v$=1, its error coefficient is zero. However, the error coefficient of CFWENO is different from FWENO and WENO, while equal to zero both at $v$=0 and $v$=1. We can also find the error coefficient of CFWENO is significantly smaller than that of FWENO and WENO.

**Table 4 Comparison of error coefficient $C_e$**

| $C_e$ | CFWENO | FWENO | WENO |
|---|---|---|---|
| 3th | $\frac{1}{(2!)^2}\lvert v\rvert(1-\lvert v\rvert)^2$ | $\frac{1}{4!}(1^2-\lvert v\rvert^2)(2-\lvert v\rvert)$ | $\frac{1}{4!}(2)$ |
| 5th | $\frac{1}{(3!)^2}\lvert v\rvert(1-\lvert v\rvert)^2(1+\lvert v\rvert)(2-\lvert v\rvert)$ | $\frac{1}{6!}(1^2-\lvert v\rvert^2)(2^2-\lvert v\rvert^2)(3-\lvert v\rvert)$ | $\frac{1}{6!}(2^2\times 3)$ |
| 7th | $\frac{1}{(4!)^2}\lvert v\rvert(1-\lvert v\rvert)^2(1+\lvert v\rvert)^2(2-\lvert v\rvert)^2$ | $\frac{1}{8!}(1^2-\lvert v\rvert^2)(2^2-\lvert v\rvert^2)(3^2-\lvert v\rvert^2)(4-\lvert v\rvert)$ | $\frac{1}{8!}(2^2\times 3^2\times 4)$ |

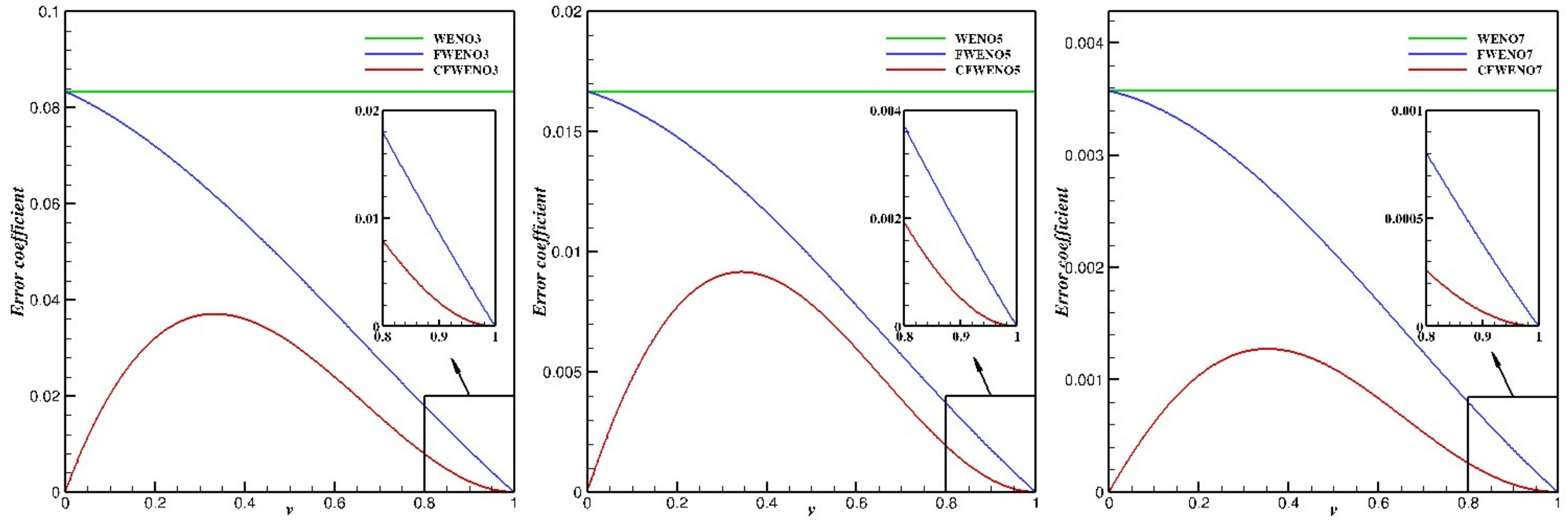

**Fig. 2 Error coefficient: comparations of 3th, 5th,7th (from left to right) order schemes**

### 3.2.3 Relevance between Hermite interpolation and compact fully-discrete WENO

Compact schemes can generally be attributed to various forms of Hermite interpolation, this is also true for CFWENO. Unlike the semi-discrete Hermite WENO (HWENO) [35–38] in Eq.(34), which uses hyperbolic conservation laws and their derivative equations for Hermite interpolation, CFWENO can be seen as using the HJ equation and hyperbolic conservation laws for Hermite interpolation. However, the final form of CFWENO does not include variables related to the HJ equation. This is because, through the solution formula method, these variables can be converted into conservative variables of the conservation laws.

$$\begin{cases} u_t + f(u)_x = 0, \quad u(x,0) = u_0(x), \\ v_t + g(u,v)_x = 0, \quad v(x,0) = v_0(x), \quad v = u_x, \\ g(u,v) = f'(u)u_x = f'(u)v. \end{cases} \tag{34}$$

The strategy of HWENO adds equations, significantly increasing the computational cost, which is even more pronounced in multi-dimensional cases. As shown in Eq.(35) for the two-dimensional case, it requires two derivative equations and also includes mixed derivatives that need to be handled separately.

$$\begin{cases} u_t + f(u)_x + g(u)_y = 0, \\ u(x,y,0) = u_0(x,y), \end{cases} \quad \begin{cases} (u_x)_t + f(u)_{xx} + g(u)_{xy} = 0, \\ (u_y)_t + f(u)_{xy} + g(u)_{yy} = 0. \end{cases} \tag{35}$$

While CFWENO does not require the computation of additional equations, making its implementation more convenient and introducing less computational cost, at least for one-dimensional cases. For multi-dimensional cases, refer to the analysis in section 3.3 below.

Additionally, the maximum accuracy of the stencils that can be constructed by CFWENO and HWENO are also different. Considering the dependency domain [$x_{i-1}$, $x_{i+1}$], CFWENO has five information, allowing for the construction of fifth-order stencil, while HWENO has six information, allowing for sixth-order stencil. However, Considering the dependency domain [$x_{i-3/2}$, $x_{i+3/2}$], CFWENO has seven information, allowing for the construction of seventh-order stencil.

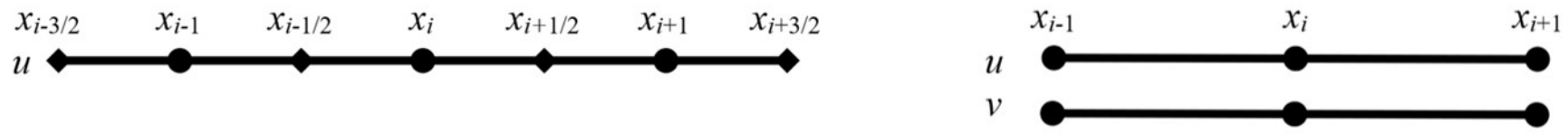

**Fig. 3 Stencils of reconstruction schemes: CFWENO (left), HWENO (right)**

## 3.3 Multi-dimensional case

For multi-dimensional cases, general difference schemes can use the dimension-by-dimension approach [57]. By this way, apart from the increase in the number of equations, the form of scheme is entirely similar to that of one-dimensional problems. The solution in any direction can be given by

$$\mathbf{u}_{\mathrm{j}}^{n+1} = \mathbf{u}_{\mathrm{j}}^{n} - \frac{\tau}{h}(\tilde{\mathbf{f}}^{e}_{\mathrm{j}+\frac{1}{2}} - \tilde{\mathbf{f}}^{e}_{\mathrm{j}-\frac{1}{2}}) \tag{36}$$

Rewrite it into following operator form

$$\mathbf{u}^{n+1} = \mathbf{T}_e(\mathbf{u}^n) \tag{37}$$

Here, the subscript $e$=$x$, $y$, $z$ indicates the direction of temporal evolution. For curvilinear grids, $e$=$\xi$, $\eta$, $\zeta$ represents the three directions of the curvilinear grid. Thus, the complete multi-dimensional solution process can be written as

$$\mathbf{u}^{n+1} = \mathbf{T}_x\mathbf{T}_y\mathbf{T}_z(\mathbf{u}^n) \tag{38}$$

For FWENO, the conventional approach can be used for temporal evolution, same as WENO+RK. However, CFWENO faces issues. As shown in Fig. 4, the dots represent grid points $u_i^n$, obtained from Eq.(14), while the diamonds represent half points $u_{i+1/2,j}^n$ or $u_{i,j+1/2}^n$ obtained from Eq.(13). The blue arrows indicate the evolution in $x$-direction, and the green arrows indicate the evolution in $y$-direction. If this traditional temporal evolution strategy is used, CFWENO will only evolve in either $x$-direction or $y$-direction at the half points, which leads to a loss of accuracy. Therefore, we designed a new evolution approach for CFWENO.

Fig. 5 shows the newly designed temporal evolution method for CFWENO, primarily aimed at repairing the accuracy loss issues. Taking $x$-direction as an example, as shown in Fig. 5, when evolving in the $x$-direction, it is necessary to calculate all the points that may be used for $y$-direction. Similarly, the dots are obtained from Eq.(14), and the diamonds are obtained from Eq.(13). For ease of programming implementation, both the dots and diamonds can be stored on the same grid. The same strategy applies to the $y$-direction evolution.

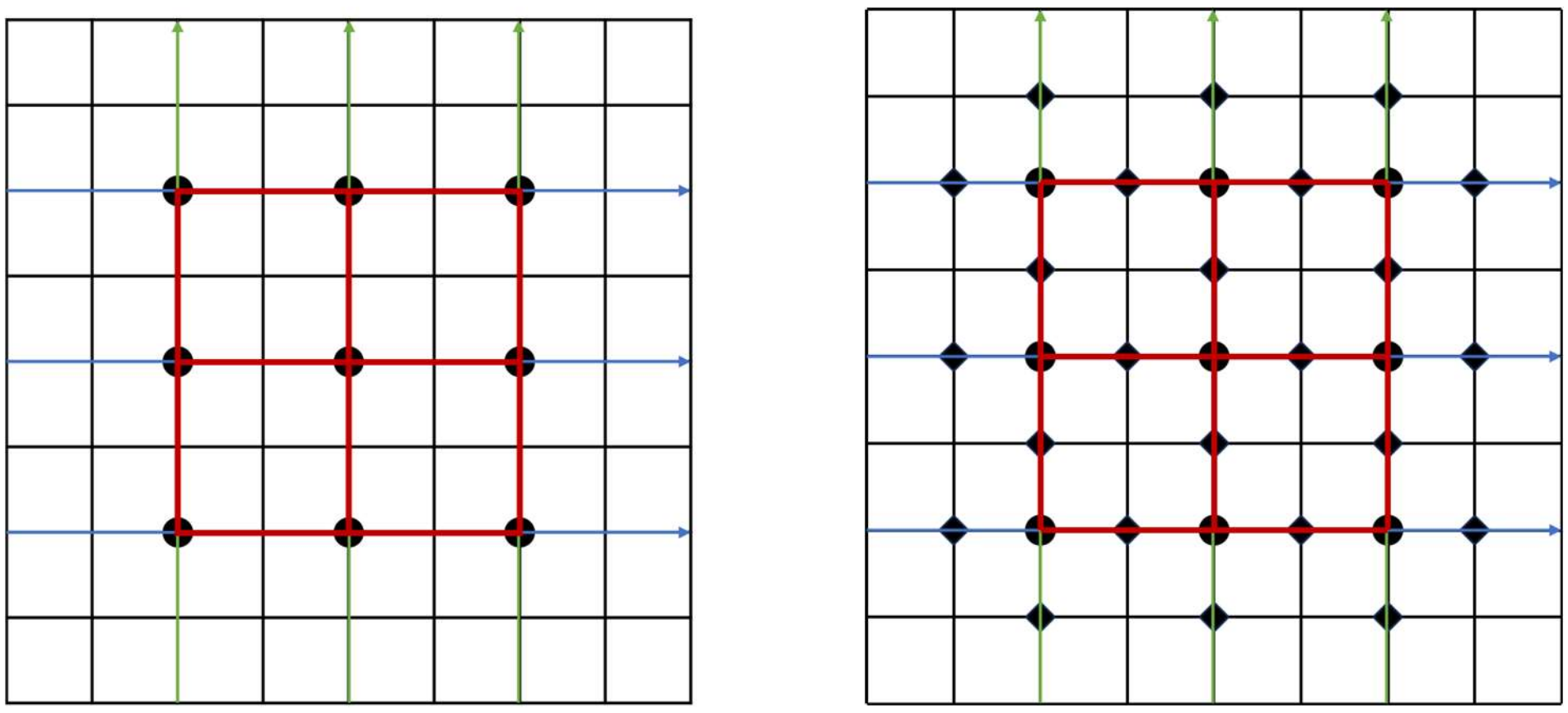

**Fig. 4 Traditional multi-dimensional strategy for FWENO/WENO (left) and CFWENO (right)**

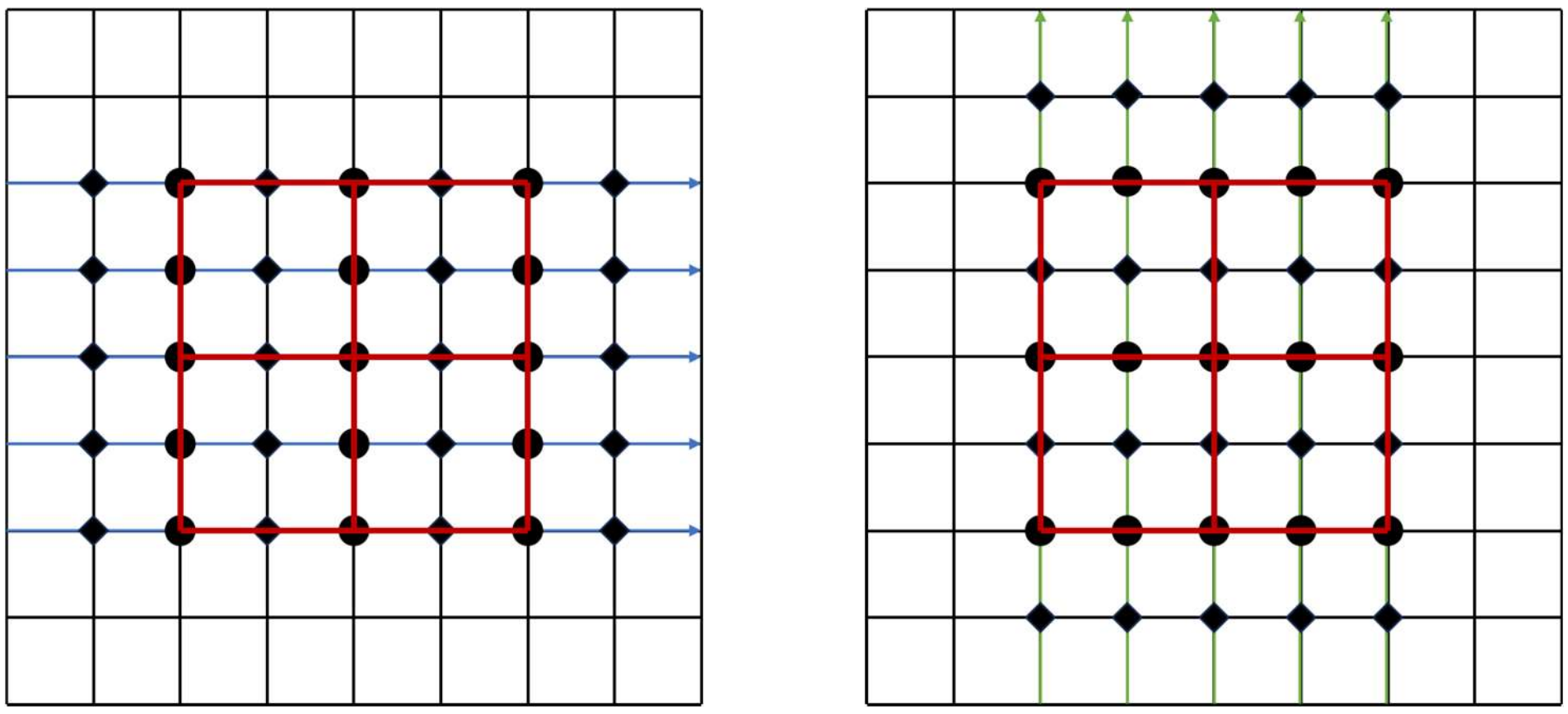

**Fig. 5 New multi-dimensional strategy for CFWENO, *x*-dimension (left) and *y*-dimension (right)**

For the new time evolution method, the computational load seems to be significantly increased. However, in the solution formula framework, the diamonds are byproducts of solving the dots, so the additional computational load is not excessive. The actual computational load only depends on the number of dots. If the computational conditions and spatial steps are the same, and if the number of grid points in each direction for FWENO/WENO is $n$, the computational load for different dimensions can be approximately given as follows:

**Table 5 Scheme computing burden in multidimensional case with same computational condition**

| $C_Q$ | CFWENO | FWENO | WENO+RK3 |
|---|---|---|---|
| 1D | $n \cdot C_{Q,\mathrm{CFWENO}}$ | $n \cdot C_{Q,\mathrm{FWENO}}$ | $n \cdot C_{Q,\mathrm{WENO}} \cdot P$ |
| 2D | $n(2n\text{-}1) \cdot C_{Q,\mathrm{CFWENO}}$ | $n^2 \cdot C_{Q,\mathrm{FWENO}}$ | $n^2 \cdot C_{Q,\mathrm{WENO}} \cdot P$ |
| 3D | $n(2n\text{-}1)^2 \cdot C_{Q,\mathrm{CFWENO}}$ | $n^3 \cdot C_{Q,\mathrm{FWENO}}$ | $n^3 \cdot C_{Q,\mathrm{WENO}} \cdot P$ |

Where $C_Q$ represents the computational cost of scheme itself obtained by numerical experiment, $P$ is the stage of temporal evolution given in Table 6. With same spatial step length and computing condition, and sufficient large grid size, the computational cost of CFWENO, FWENO and WENO+RK3 are given in Table 5, which shows: (1) For 1D case, CFWENO is $C_{Q,\mathrm{CFWENO}}$ /$C_{Q,\mathrm{FWENO}}$ times of FWENO, and $C_{Q,\mathrm{CFWENO}}/(C_{Q,\mathrm{WENO}} \cdot P)$ times of WENO+RK3; (2) For 2D case, CFWENO is $2C_{Q,\mathrm{CFWENO}}/C_{Q,\mathrm{FWENO}}$ of FWENO, and $2C_{Q,\mathrm{CFWENO}}/(C_{Q,\mathrm{WENO}} \cdot P)$ of WENO+RK3; (3) For 3D case, CFWENO is $4C_{Q,\mathrm{CFWENO}}/C_{Q,\mathrm{FWENO}}$ of FWENO, and $4C_{Q,\mathrm{CFWENO}}/(C_{Q,\mathrm{WENO}} \cdot P)$ of WENO+RK3

The comparison above considers the computational cost under same spatial step length in the evolution direction. Another comparison standard is to regard both the dots and diamonds as grid points to compare computational efficiency. Then, under same grid resolutioin, the computing speed can be given by following formula:

$$Q_e = \frac{\Delta_e \bullet \Delta_{t,e}}{C_Q \bullet P}, \quad \Delta_{t,e} = \frac{\Delta_e \bullet \text{CFL}}{a_{e,\max}} \tag{39}$$

Where the subscript $e$ represents for evolution direction, $Q_e$ is computing speed, $C_Q$ is the computational cost of scheme, $P$ is the stage of evolution, $\Delta_e$ is the time step length, $a_{e,max}$ is the global maximum eigenvalue which can be seen as the same in same case. The CFL number for CFWENO and FWENO schemes is 0.9, while it is commonly 0.6 for WENO + RK3. The specific parameters are shown in Table 6.

$C_Q$ shown in Table 7 is an average value calculated by two 1D Euler numerical test. With these, we can estimate the computing speed of three schemes. Table 8 gives the estimated computational speed $Q_e$, which shows the computing speed of CFWENO reaches about three times of FWENO, and even higher than ten times of WENO+RK3, if counts both dots and diamonds into grid points. This also has been verified in section 4.5.

**Table 6 Parameters for computational speed $Q_e$**

| | CFWENO | FWENO | WENO+RK3 |
|---|---|---|---|
| $P$ | 1 | 1 | 3 |
| $\Delta_e$ | 2 | 1 | 1 |
| CFL | 0.9 | 0.9 | 0.6 |

**Table 7 Scheme computing burden $C_Q$ for computational speed $Q_e$**

| $C_Q$ | CFWENO | FWENO | WENO+RK3 |
|---|---|---|---|
| 3th | 1.71 | 1.29 | 1 |
| 5th | 1.47 | 1.10 | 1 |
| 7th | 1.31 | 0.97 | 1 |

**Table 8 Computational speed $Q_e$ (normalized $Q_e$)**

| $Q_e$ | CFWENO | FWENO | WENO+RK3 |
|---|---|---|---|
| 3th | 2.11 (10.53) | 0.70 (3.49) | 0.20 (1.00) |
| 5th | 2.45 (12.24) | 0.82 (4.09) | 0.20 (1.00) |
| 7th | 2.75 (13.74) | 0.93 (4.64) | 0.20 (1.00) |

# 4 Numerical verifications

This section considers various benchmark examples to test and evaluate the new scheme, mainly including accuracy tests, scalar equations, one-dimensional and multi-dimensional Euler equations. In this paper, we set FWENO and WENO+RK3 as a comparison. If not specified, CFWENO and FWENO use CFL=0.9, and WENO+RK3 uses CFL=0.6, while FWENO also uses the entropy condition flux reconstruction method design in Eq.(30). In addition, due to CFWENO and FWENO applies Roe average for discontinuities in flux reconstruction, we also use Roe flux with entropy fix for WENO+RK3.

In following tests, the grid points of scalar equation and one-dimensional Euler equations mean the quantity of node value $u_i$, which does not take the half point $u_{i+1/2}$ into account. For multi-dimensional cases, due to the newly designed evolution method, the mesh resolution represents the sum of all node points and half points.

## 4.1 Accuracy test

This section compares the accuracy and errors of CFWENO, FWENO and WENO+RK.

**Example 4.1.1 Linear equation**

Consider following linear scalar equation

$$\begin{cases} u_t + u_x = 0, \\ u_0(x) = \sin(\pi x), \end{cases} \quad x \in [-1,1], \quad t = 2 \tag{40}$$

Fig. 6 shows CFWENO and FWENO can basically satisfy designed accuracy for linear flux, while semi-discrete WENO only can reach designed accuracy when using RK method with corresponding accuracy. In general, to eliminate the accuracy impact of RK method, a very small CFL generally need to be implemented with low-order RK method for accuracy test. Moreover, compared to FWENO and WENO+RK, although they all use the WENO reconstruction, CFWENO significantly has lower errors due to its compact nature.

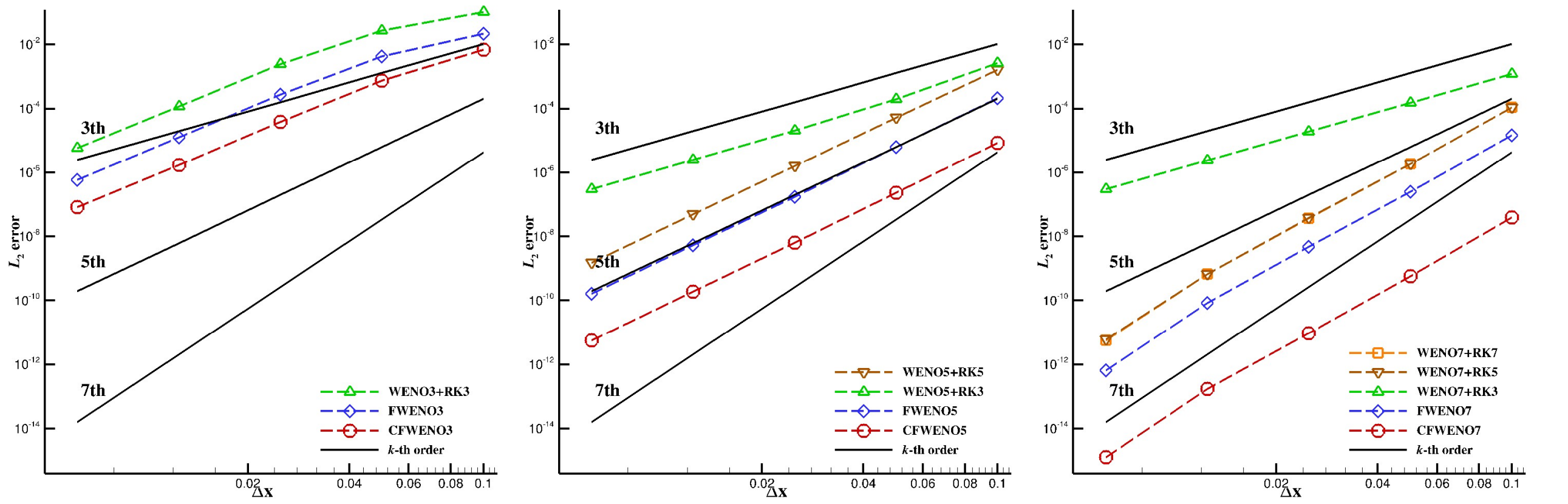

**Fig. 6 Convergence of the $L_2$ error for linear equation (CFL=0.5): comparations of 3th, 5th,7th (from left to right) order schemes**

**Example 4.1.2 Nonlinear equation**

Consider following nonlinear scalar equation

$$\begin{cases} u_t + (\frac{u^2}{2})_x = 0, \\ u_0(x) = 0.5 + \sin(\pi x), \end{cases} \quad x \in [0,2], \quad t = 0.15 \tag{41}$$

Fig. 7 shows that, similar to the linear case, both CFWENO and FWENO essentially achieve the designed accuracy. In contrast, the semi-discrete WENO requires the corresponding accuracy RK method to achieve designed accuracy. Moreover, the errors of CFWENO at each order are lower than those of FWENO and WENO+RK.

Fig. 8 and Fig. 9 show error variation with respect to the number of flux reconstruction iterations for FWENO and CFWENO. In Fig. 8 and Fig. 9, iteration=0 indicates directly using baseline entropy condition flux reconstruction, iteration=1 indicates using the eigenvalues from iteration=0 to obtain the new eigenvalue, iteration=2 indicates using the eigenvalues from iteration=1 to obtain the new eigenvalue, and higher iteration counts follow this pattern. Here we can find, as analyzed previously in section 3.2.1, the error of FWENO and CFWENO decreases with the increasing iteration numbers. FWENO improves by two orders of accuracy with each iteration, while CFWENO improves by only one order.

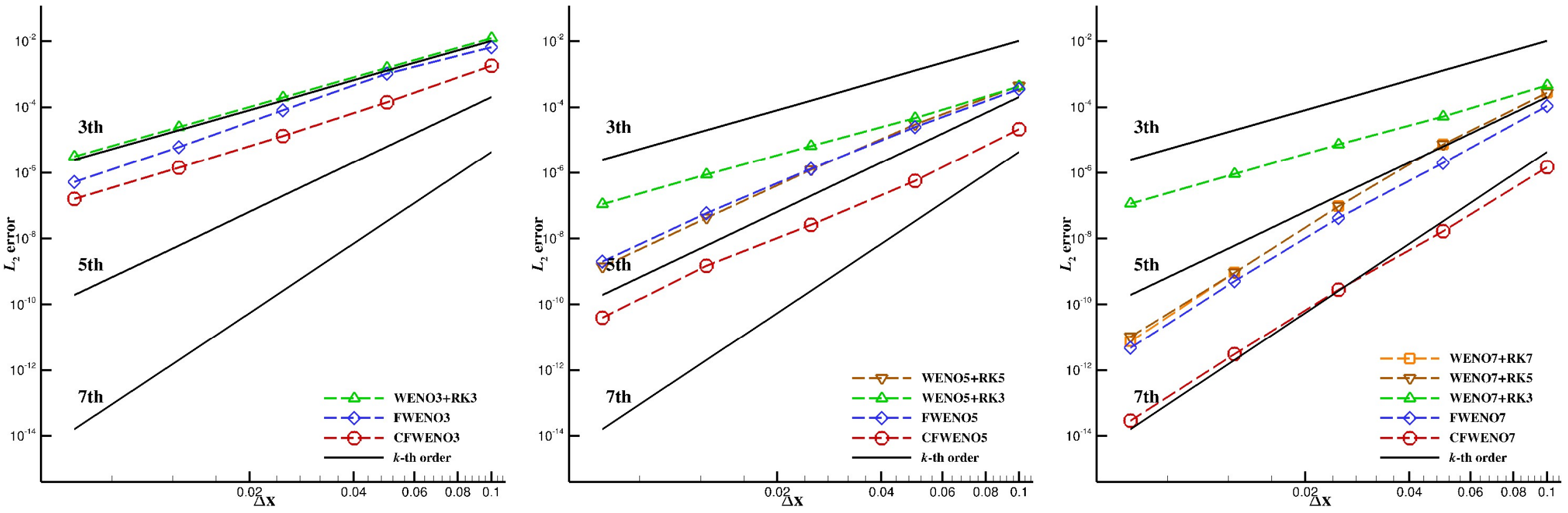

**Fig. 7 Convergence of the $L_2$ error for nonlinear equation: comparations of 3th, 5th,7th (from left to right) order schemes**

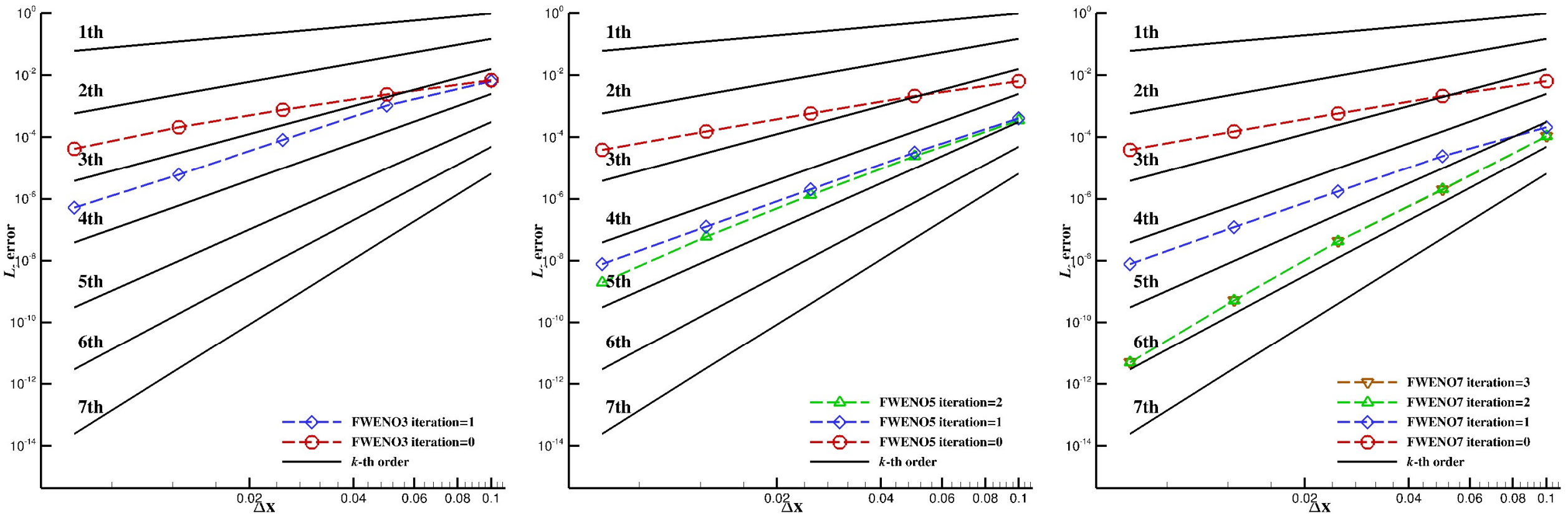

**Fig. 8 Convergence of the $L_2$ error for nonlinear equation: Comparison of 3th, 5th,7th (from left to right) order schemes**

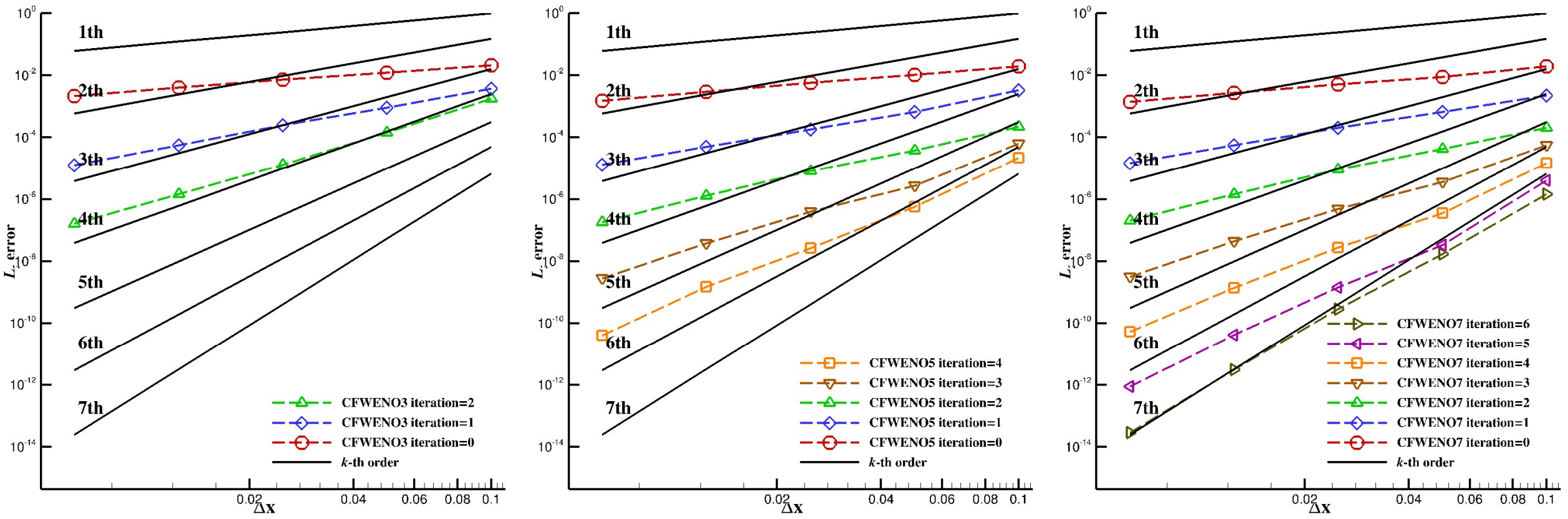

**Fig. 9 Convergence of the $L_2$ error for nonlinear equation: Comparison of 3th, 5th,7th (from left to right) order schemes**

## 4.2 Scalar equation

### Example 4.2.1 Linear scalar equation: Square wave problem

We now consider the square wave problem for linear scalar equation. This test aims to show the different properties of these three schemes with increasing CFL number. The initial condition is

$$\begin{cases} u_t + u_x = 0, \\ u_0(x) = \begin{cases} 1, & -\frac{1}{3} \le x \le \frac{1}{3}, \\ -1, & else, \end{cases} \quad x \in [-1,1]. \end{cases} \tag{42}$$

The simulation is performed on a uniform mesh with $N$=100 and output at simulation time $t$=20. It can be shown in Fig. 10, the resolution of FWENO rises along with the increasing CFL number, while WENO+RK3 shows the opposite behavior. As for CFWENO, it also shows rising resolution with increasing CFL number, although the error coefficient monotonically decreasing. This is because there are more time steps when using a small CFL number, which would introduce cumulative error. Moreover, CFWENO and FWENO obtain an exact solution when CFL=1, because the solution formula method is exact for linear flux.

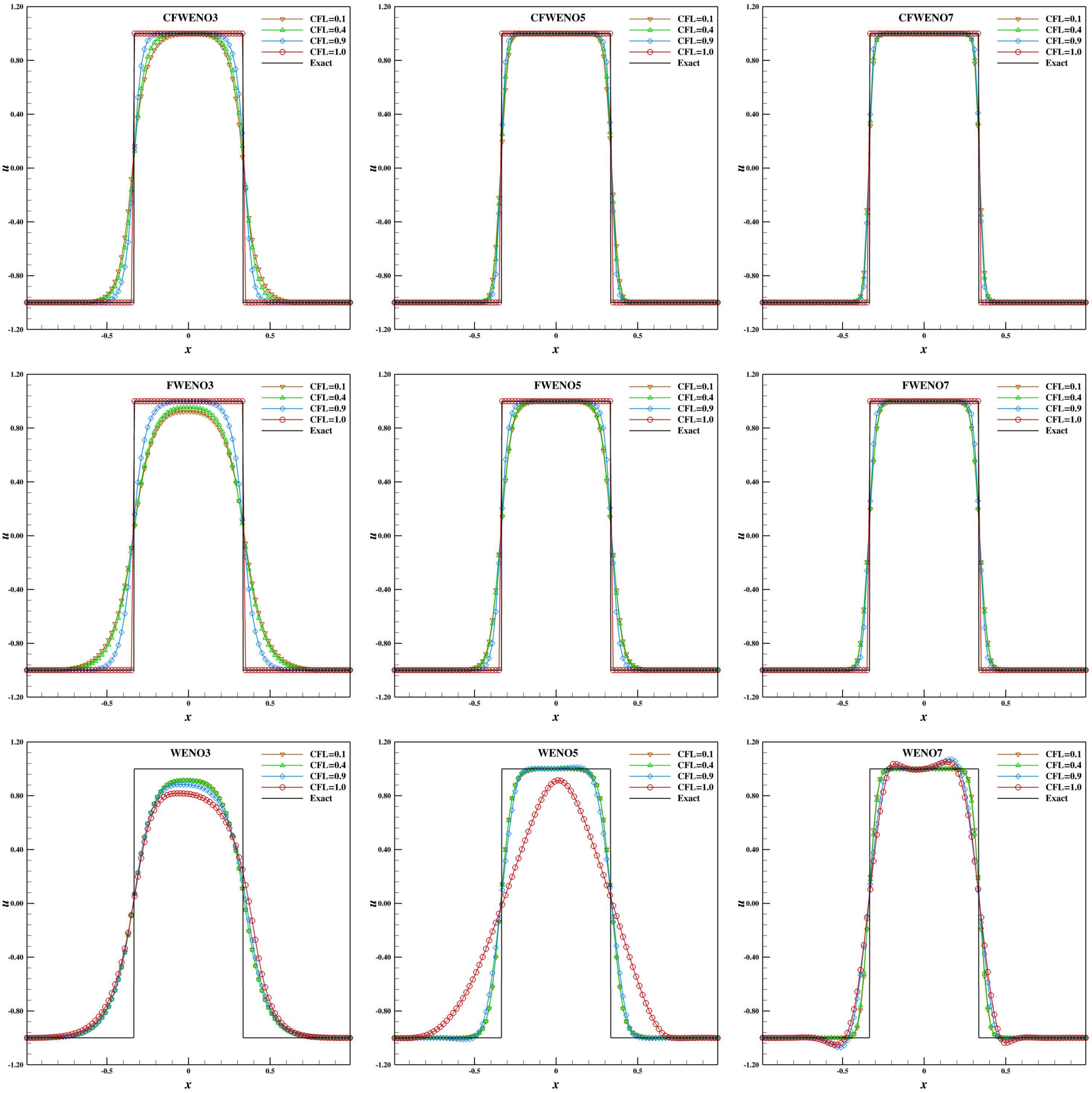

**Fig. 10 Square wave problem: solution at simulation t=20. Comparison of CFWENO, FWENO, WENO+RK3 (from top to below) and 3th, 5th, 7th order schemes (from left to right). Discretization on 100 uniformly distributed grid points.**

**Example 4.2.2 Linear scalar equation: Multiple extremes problem [7]**

This test involves a combination of smooth and non-smooth functions, including a Gaussian, a square, a triangle, and a semi-ellipse. So, it is widely used to test the performance of discontinuity capturing. The initial condition is

$$\begin{cases} u_t + u_x = 0, \\ u_0(x) = \begin{cases} \frac{1}{6}(G(x,\beta,z-\delta) + G(x,\beta,z+\delta) + 4G(x,\beta,z)), & \text{-}0.8 \le x \le \text{-}0.6, \\ 1, & \text{-}0.4 \le x \le \text{-}0.2, \\ 1 - |10(x-0.1)|, & 0 \le x \le 0.2, \\ \frac{1}{6}(F(x,\alpha,a-\delta) + F(x,\alpha,a+\delta) + 4F(x,\alpha,a)), & 0.4 \le x \le 0.6, \\ 0, & \text{else}, \end{cases} \\ G = e^{-\beta(x-z)^2}, \quad F(x,\alpha,a) = \sqrt{\max(1-\alpha^2(x-a)^2, 0)}, \\ z = -0.7, \quad \delta = 0.005, \quad \beta = (\log 2)/(36\delta^2), \quad a = 0.5, \quad \alpha = 10. \end{cases} \tag{43}$$

The simulation is performed on a uniform mesh with $N$=200 and output at simulation time $t$=8, 80, 800 respectively. In Fig. 11, it can be observed that the resolution of CFWENO is significantly higher than that of FWENO and WENO+RK3. This phenomenon becomes more pronounced over time: at $t$=8, all three methods can still resolve the waves well; by $t$=80, due to the dissipation along with time evolution, the resolution of WENO+RK3 noticeably decreases, and FWENO also shows significantly lower resolution than that of CFWENO at extremes; by $t$=800, WENO+RK3 shows obvious distorted, and resolution of FWENO also greatly decreases, but benefit from its compact property and extremely small error coefficient, CFWENO can still resolve the waves well.

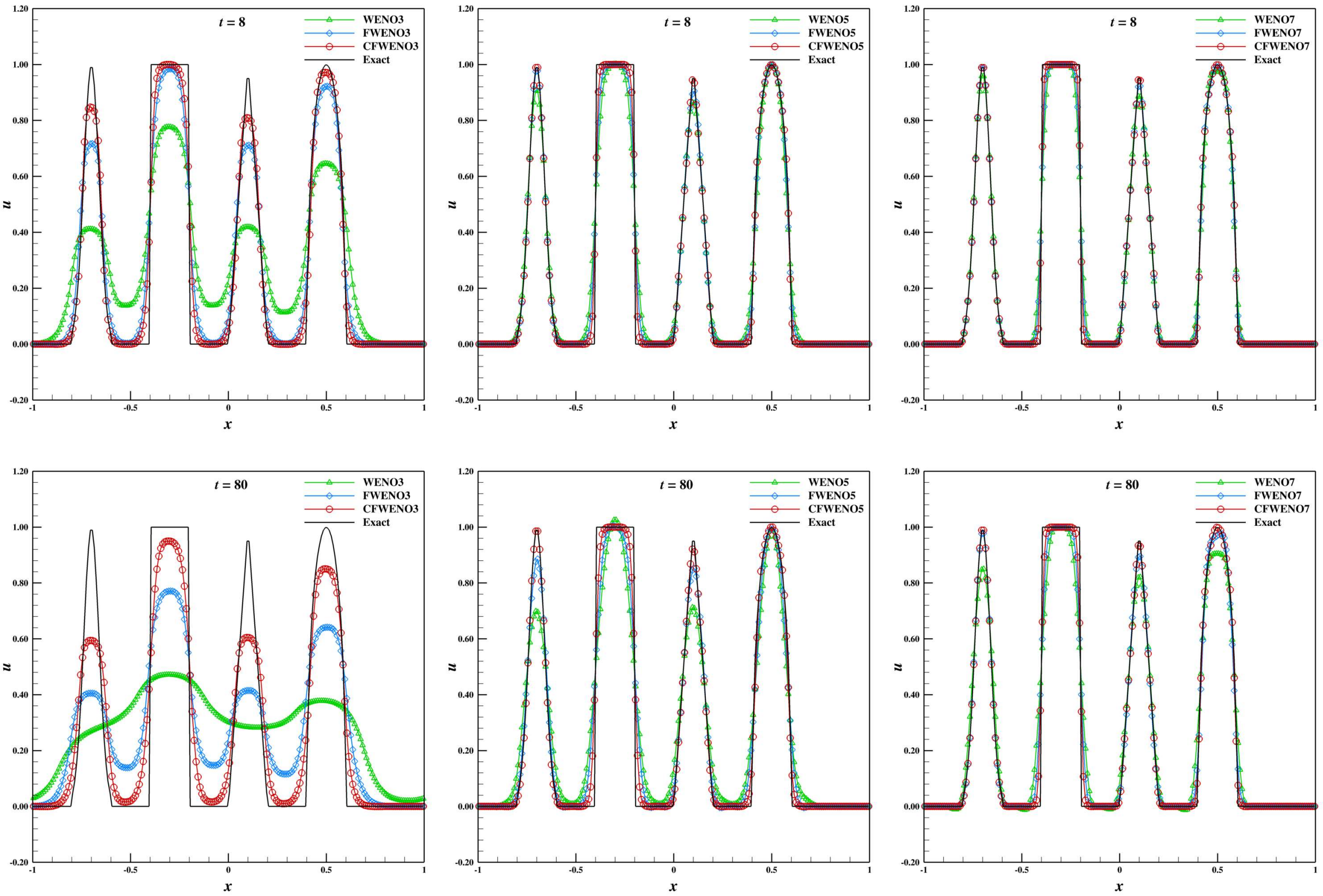

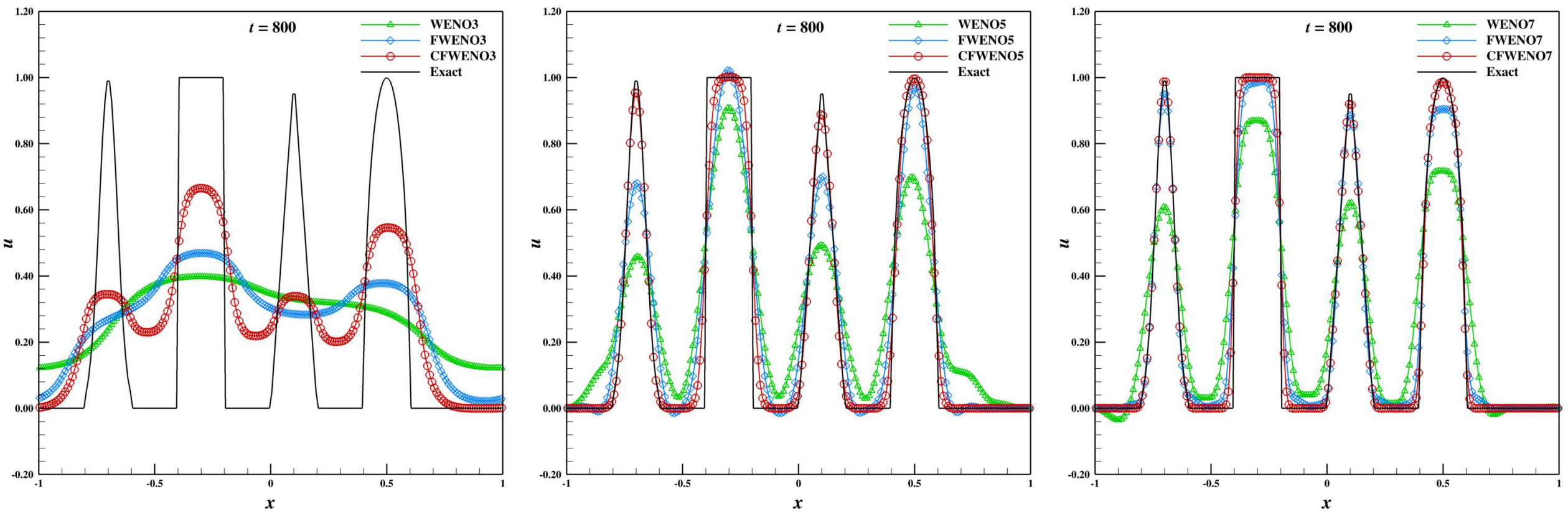

**Fig. 11 Multiple extremes problem: solution at simulation t=8, 80, 800 (from top to below). Comparison of 3th, 5th, 7th order schemes (from left to right). Discretization on 200 uniformly distributed grid points.**

**Example 4.2.3 Nonlinear scalar equation: Burgers equation**

Here we consider the nonlinear Burgers equation Eq.(41). The simulation is performed on a uniform mesh with $N$=80 and output at simulation time $t$=20, 200, 2000 respectively, here CFWENO and FWENO use CFL=1. Even if the initial values of the Burgers equation are sufficiently smooth, discontinuities will appear when time evolves. In Fig. 12, it can be observed that at $t$=20, all three methods can solve the discontinuity well, while CFWENO behaves best, followed by FWENO, and finally WENO+RK3. By $t$=200, WENO+RK3 occurs suspicious oscillations near the discontinuity, but CFWENO and FWENO still maintain high resolution. By $t$=2000, both FWENO and WENO+RK3 generate suspicious oscillations, but CFWENO still maintains very high resolution at discontinuity, which indicates that CFWENO has lower dissipation over time.

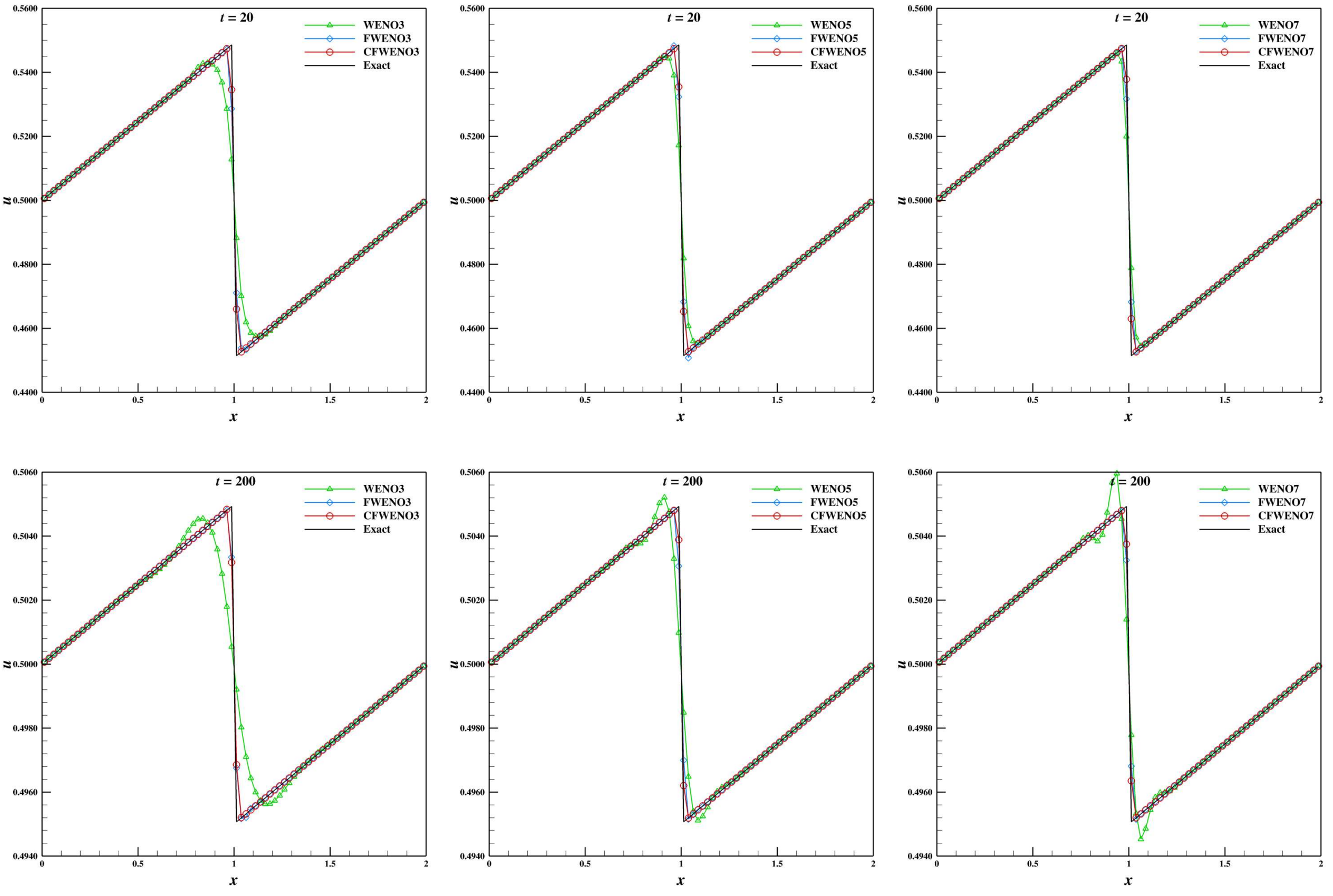

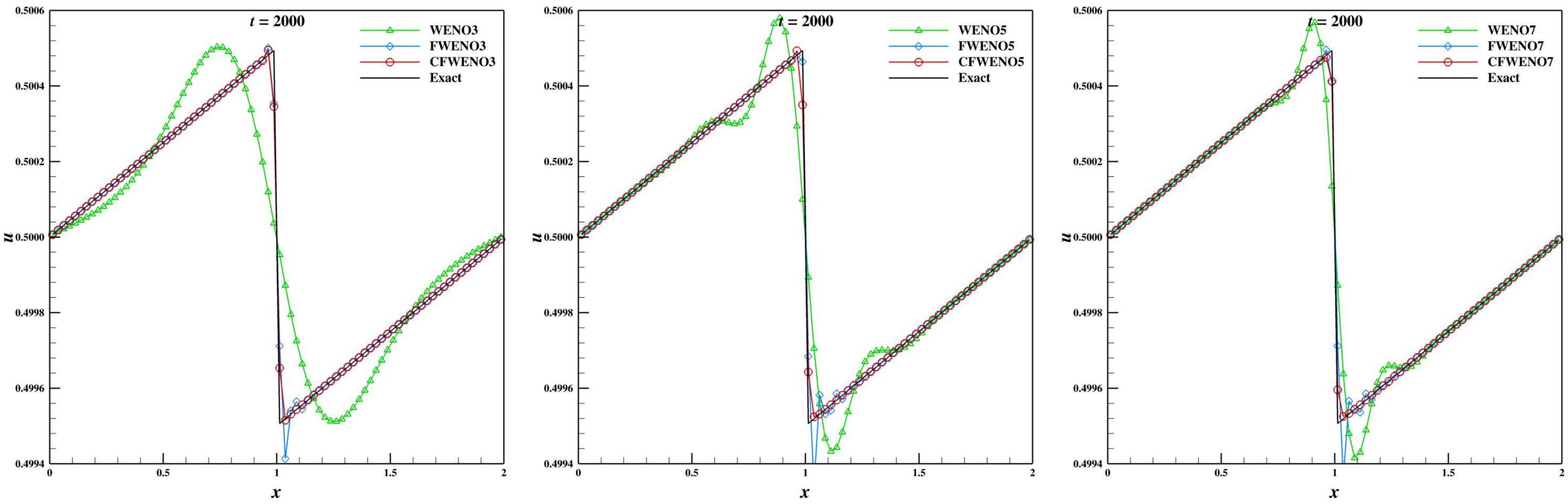

**Fig. 12 Burgers equation: solution at simulation time t=20, 200, 2000 (from top to below). Comparison of 3th, 5th, 7th order schemes (from left to right). Discretization on 80 uniformly distributed grid points.**

## 4.3 One-dimensional Euler equation

**Example 4.3.1 Sod shock tube problem** [59]

The typical Riemann problem Sod shock tube problem is considered here. The initial conditions discontinued in the interface would evolve into a smooth expansion wave, a discontinuous entropy wave (contact discontinuity) and a discontinuous shock wave. Whether the discontinuities can be captured correctly, particularly for the contact discontinuities in density profile, is an important index to measure the discontinuity capturing ability of a scheme. The initial condition is

$$(\rho, u, p)=\begin{cases}(1, 0, 1), & 0\leq x\leq 0.5,\\ (0.125, 0, 0.1), & 0.5\leq x\leq 1.\end{cases} \tag{44}$$

The simulation is performed on a uniform mesh with $N$=200 and output at simulation time $t$=0.2. The computed density profiles are shown in Fig. 13. None of these three schemes exhibit significant oscillations near the discontinuities, but CFWENO resolves better at both contact discontinuities and shock waves, with fewer transition points. As shown in section 4.5, Table 9, Table 10, and Table 11 reveal that the computational cost of CFWENO is approximately 1.2 to 1.4 times that of FWENO, and about 3/10 to 4/10 of WENO+RK3. This also confirms that, in the one-dimensional case, the complexity of CFWENO compared to FWENO has not increased significantly; it only stores variables that FWENO would have needed to compute and uses them for interpolation.

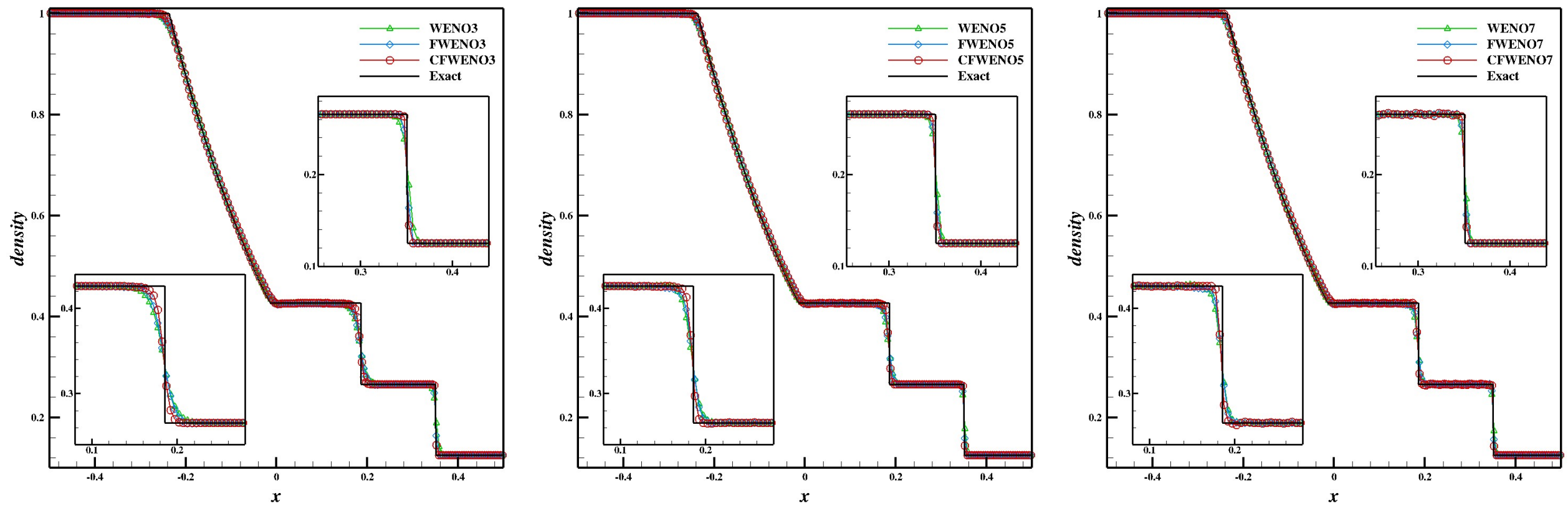

**Fig. 13 Sod problem: density distribution at simulation time *t* = 0.2. Comparison of 3th, 5th, 7th order schemes (from left to right). Discretization on 200 uniformly distributed grid points**

**Example 4.3.1 Shu-Osher problem [2]**

This test describes the interaction between a one-dimensional Mach 3 shock wave and a perturbed density field. The interaction will generate both small-scale structures and discontinuities, hence is often selected to validate shock-capturing and wave-resolution capability. The initial condition is

$$(\rho, u, p)=\begin{cases}(3.857, 2.629, 10.333), & -5\le x<-4,\\ (1+0.2\sin(5x), 0, 1), & -4\le x<5.\end{cases} \tag{45}$$

The simulation is performed on a uniform mesh with $N$=200 and output at simulation time $t$=1.8. The reference solution is computed by the fifth-order WENO5-JS with $N$=20000. As shown in Fig. 14, for different order schemes, FWENO3 and WENO+RK3 have comparable resolution, with CFWENO being significantly higher than both, especially at the extremes.

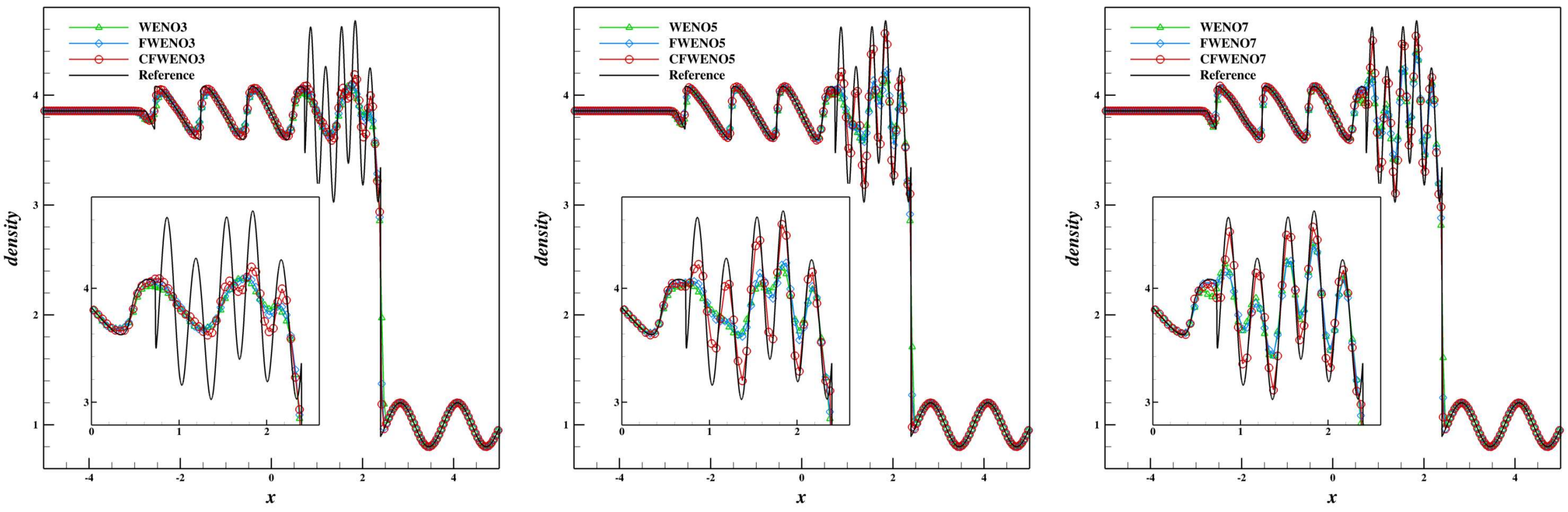

**Fig. 14 Shu-Osher problem: density distribution at simulation time $t$ = 1.8. Comparison of 3th, 5th, 7th order schemes (from left to right). Discretization on 200 uniformly distributed grid points**

**Example 4.3.2 Titarev-Toro problem** [60]

This is an extending Shu-Osher problem designed by Titarev and Toro [60], which is always used to examine the performance at the interaction of high-frequency waves and shocks. The initial condition is

$$(\rho, u, p)=\begin{cases}(1.515695, 0.523346, 1.80500), & -5\le x\le -4.5,\\ (1+0.1\sin 20\pi x, 0, 1), & -4.5\le x\le 5.\end{cases} \tag{46}$$

The simulation is performed on a uniform mesh with both $N$=400 and $N$=800. We output the solution at simulation time $t$=5. The reference solution is computed by the fifth-order WENO5-JS with $N$=20000.

As shown in Fig. 15, when $N$=400, FWENO and WENO+RK3 have similar resolution at extremes. However, for both the high-frequency waves with discontinuities at $x$<2 and the high-amplitude and high-frequency smooth waves at 2<$x$<3.2, CFWENO performs better. At $x$<2, FWENO5 and WENO5+RK3 cannot resolve the waves at all, while FWENO7 and WENO7+RK3 can only resolve a little. CFWENO5 and CFWENO7, however, show much higher resolution. When $N$=800, the resolution of FWENO and WENO+RK3 does not significantly improve compared to $N$=400, especially in the 2<$x$<3.2 region, while CFWENO improves significantly. Although CFWENO3 still cannot resolve well in the 2<$x$<3.2 region due to its lower accuracy, it can now discern the waves in the $x$<2 region, which FWENO3 and WENO3+RK3 cannot achieve.

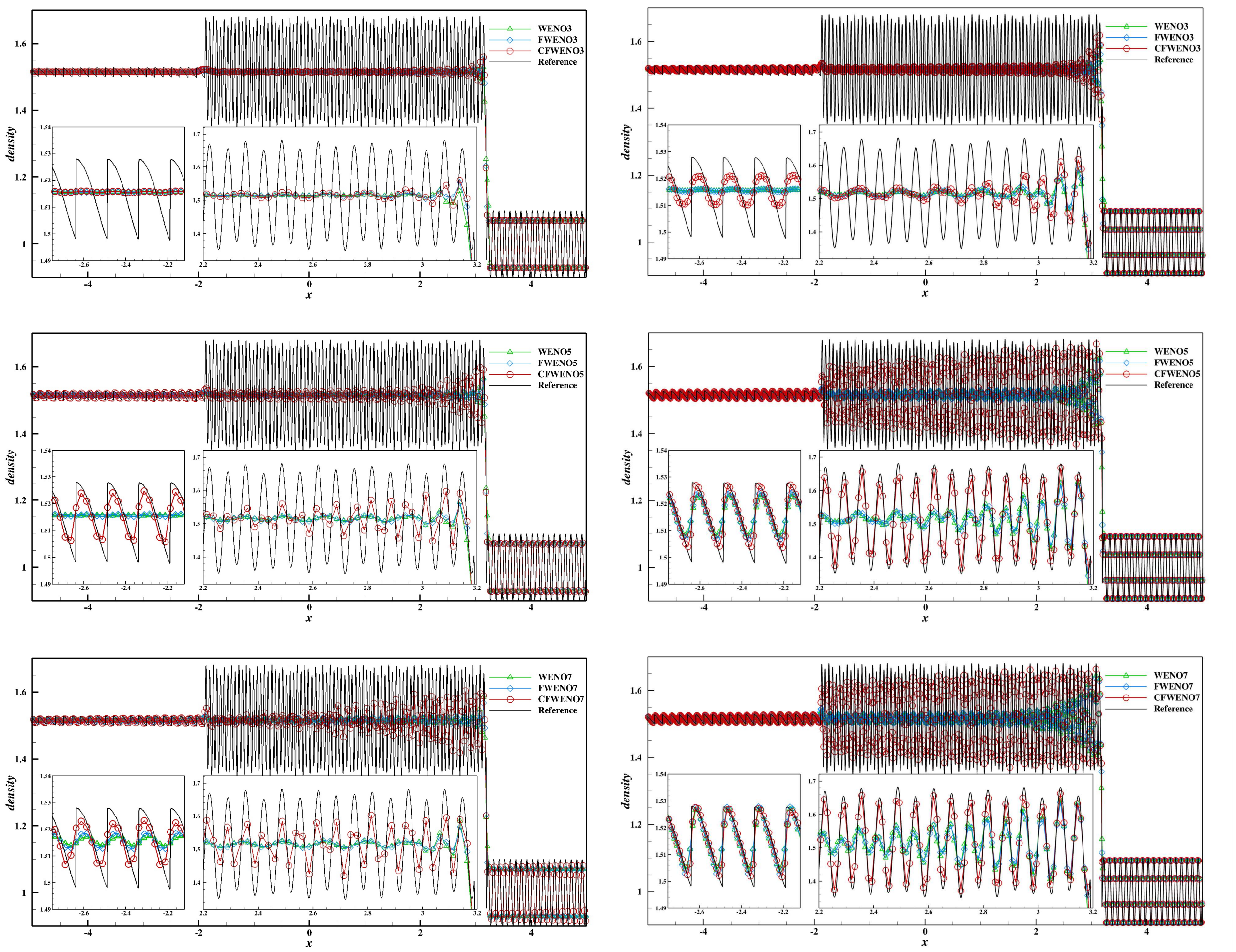

**Fig. 15 Titarev-Toro problem: density distribution at simulation time $t = 5.0$. Comparison of 3th, 5th, 7th order schemes (from top to below). Discretization on 400 (left) and 800 (right) uniformly distributed grid points**

**Example 4.3.2 Blast-wave problem** [61]

This example describes a blast-wave problem constructed by Woodward and Colella [61], which is a very severe case with very strong discontinuities for high order schemes. Reflective boundary condition is applied at $x$=0 and $x$=1, and the initial condition is

$$(\rho, u, p)=\begin{cases}(1, 0, 10^{3}), & 0\le x<0.1,\\ (1, 0, 10^{-2}), & 0.1\le x<0.9,\\ (1, 0, 10^{2}), & 0.9\le x<1.\end{cases} \tag{47}$$

The simulation is performed on a uniform mesh with $N$=200 and output at simulation time $t$=0.038. The reference solution is computed by the fifth-order WENO5-JS with $N$=20000. As shown in Fig. 16, at the contact discontinuity at $x$=0.6, the valley at $x$=0.75, and the peak at $x$=0.78, FWENO has slightly higher resolution than WENO+RK3, while CFWENO is significantly better than both. Additionally, near the shock at $x$=0.85, none of the schemes exhibit noticeable oscillations or overshoots.

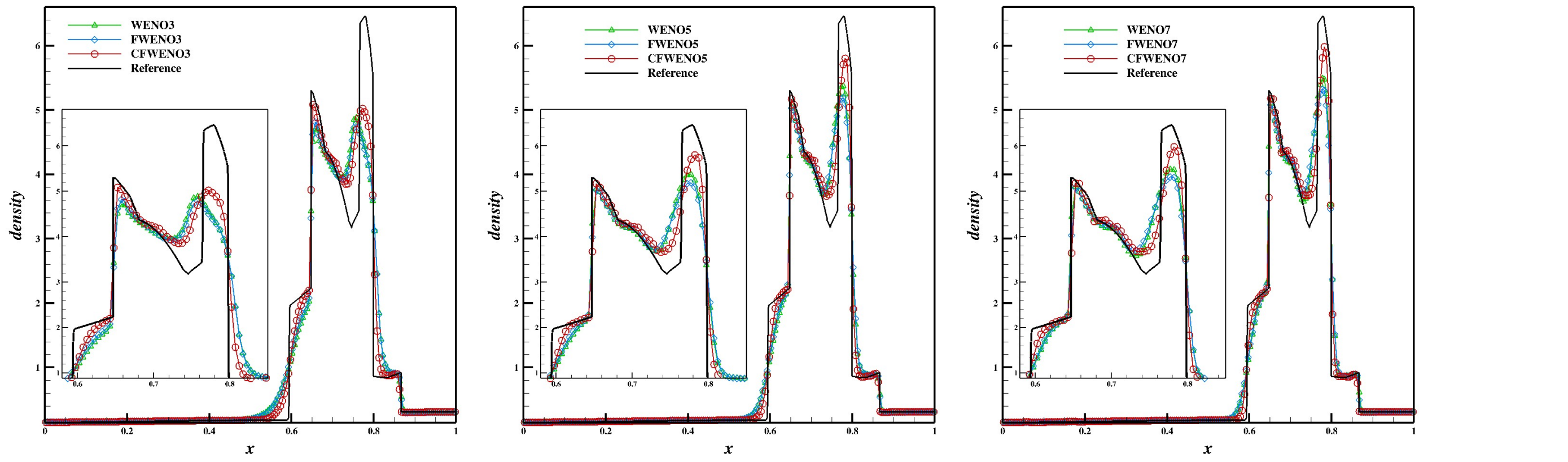

**Fig. 16 Blast-wave problem: density distribution at simulation time *t* = 0.038. Comparison of 3th, 5th, 7th order schemes (from left to right). Discretization on 200 uniformly distributed grid points**

## 4.4 Two-dimensional Euler equation

**Example 4.4.1 Two-dimensional Riemann problem 1: interaction of planar shock** [62]

We first consider the 2D Riemann problem of configuration 3 in [62], which includes interaction of planar shocks. The computational domain is set as [0,1]×[0,1] with the mesh resolution of 400×400. With the final simulation time t=0.8, the initial condition is

$$(\rho, u, v, p)=\begin{cases}(1.5, 0.0, 0.0, 1.5), & x>0.8,\ y>0.8,\\ (0.5323, 1.206, 0.0, 0.3), & x<0.8,\ y>0.8,\\ (0.138, 1.206, 1.206, 0.029), & x<0.8,\ y<0.8,\\ (0.5323, 0.0, 1.206, 0.3), & x>0.8,\ y<0.8,\end{cases} \tag{48}$$

As shown in Fig. 17, compared with FWENO and WENO + RK3, CFWENO captures more shockwave patterns and vortex structures. However, due to the small numerical viscosity of CFWENO, some numerical noise is generated, which needs to be improved in the future.

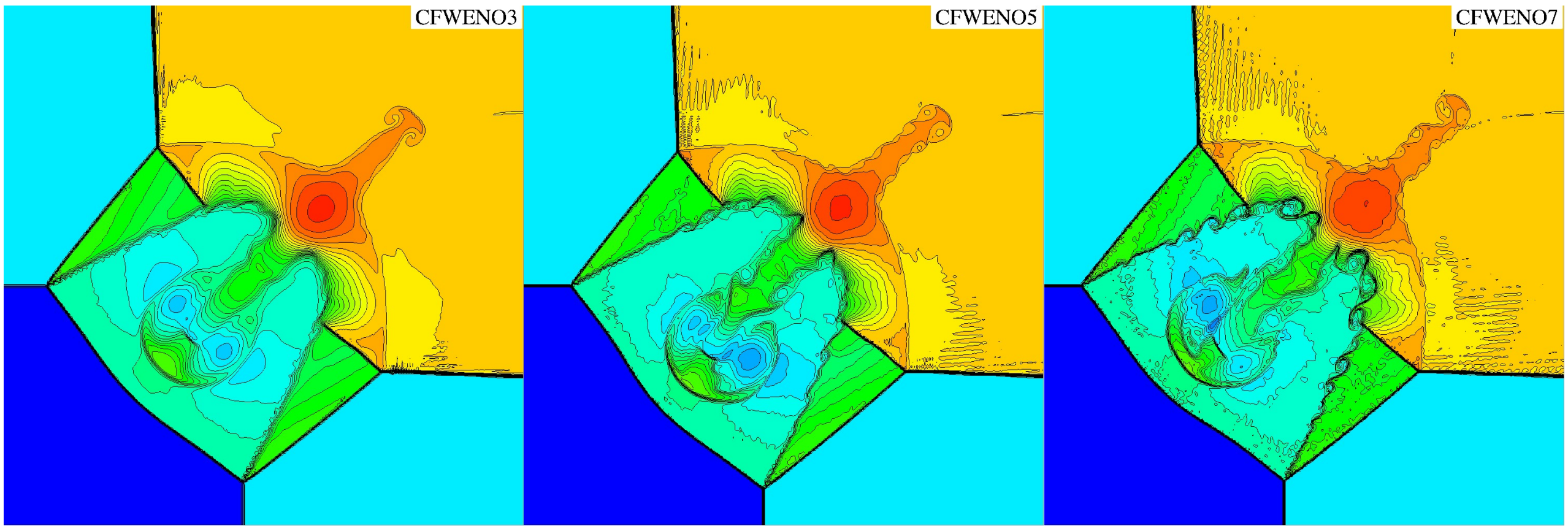

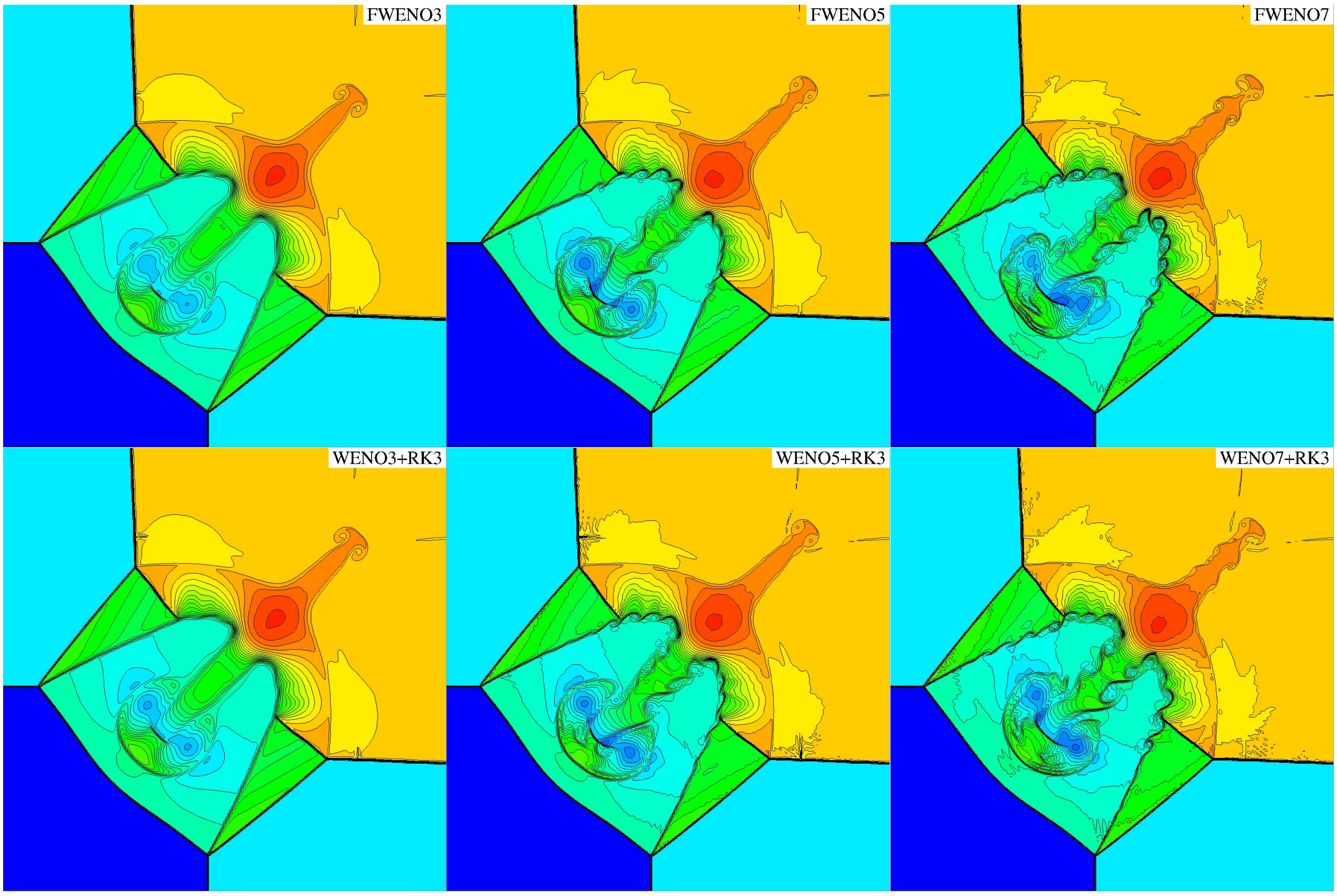

**Fig. 17 Two-dimensional Riemann problem: density contours of CFWENO, FWENO, WENO+RK3 (from top to below) and 3th, 5th, 7th order schemes (from left to right). Resolution 400×400. The simulation is run until $t$ = 0.8.**

**Example 4.4.2 Two-dimensional Riemann problem 2** [62]

We next consider the 2D Riemann problem of configuration 16 in [62]. The computational domain is set as [0,1]×[0,1] with the mesh resolution of 800×800. With finial simulation time $t$=0.6, the initial condition is

$$(\rho, u, v, p)=\begin{cases}(0.5313, 0.1, 0.1, 0.4), & x>0.5,\ y>0.5,\\(1.0222, -0.6179, 0.1, 1.0), & x<0.5,\ y>0.5,\\(0.8, 0.1, 0.1, 1.0), & x<0.5,\ y<0.5,\\(1.0, 0.1, 0.8276, 1.0), & x>0.5,\ y<0.5.\end{cases} \tag{49}$$

It is shown in Fig. 18 that CFWENO is similar to FWENO and better than WENO+RK3 in capturing the tiny structure induced by the unstable interface in the third and fifth-order cases in this example. In the seventh-order case, CFWENO captures significantly more tiny structures than the other two.

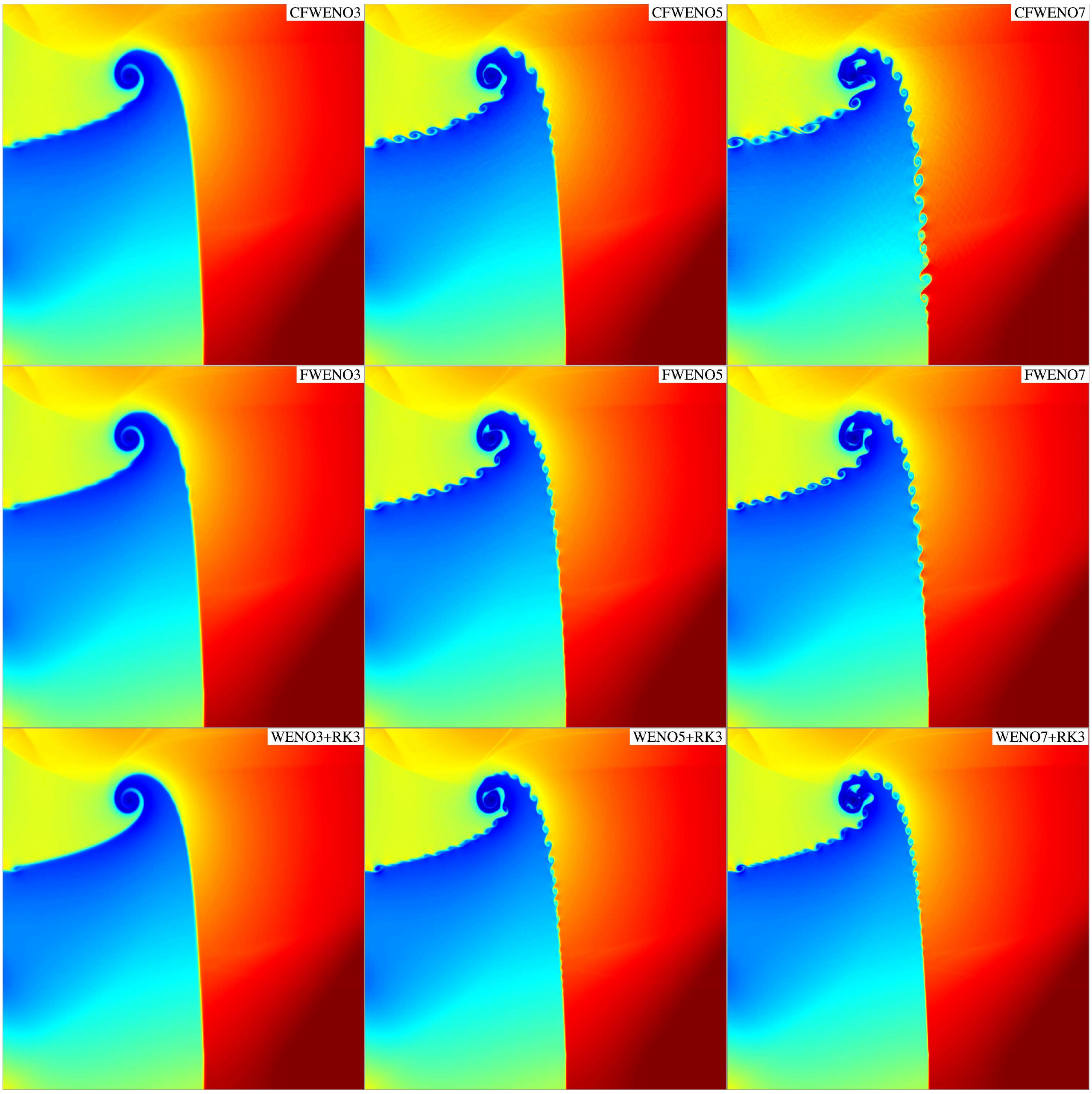

**Fig. 18 Two-dimensional Riemann problem: density contours of CFWENO, FWENO, WENO+RK3 (from top to below) and 3th, 5th, 7th order schemes (from left to right). Resolution 800×800. The simulation is run until $t$ = 0.6.**

**Example 4.4.3 Two-dimensional implosion problem** [63]

We further consider the 2D implosion problem of Liska and Wendroff [63] to test the numerical robustness of the proposed scheme. In this benchmark case, the high-density and high-pressure fluids will implode into the central region with low-density and low-pressure, generating strong shock waves propagating back and forth in the periodic box and flow instabilities at the interface between high- and low-density regions. The computational domain is set as [0, 0.6] ×[0, 0.6] with the mesh resolution of 600×600 and periodic conditions enforced for all boundaries. With the final simulation time $t$=0.4, the initial condition is

$$(\rho, u, v, p)=\begin{cases}(0.125, 0.0, 0.0, 0.14), & |x-0.3|<0.15 \text{ and } |y-0.3|<0.15, \\ (1,0,0,1), & \text{else,}\end{cases} \tag{50}$$

As shown in Fig. 19, the strong shock wave problem simulated by CFWENO has a similar structure to FWENO and WENO+RK. In addition, compared with the same-order scheme, CFWENO captures the most fine structures caused by interface instability, followed by FWENO, and finally WENO+RK3.

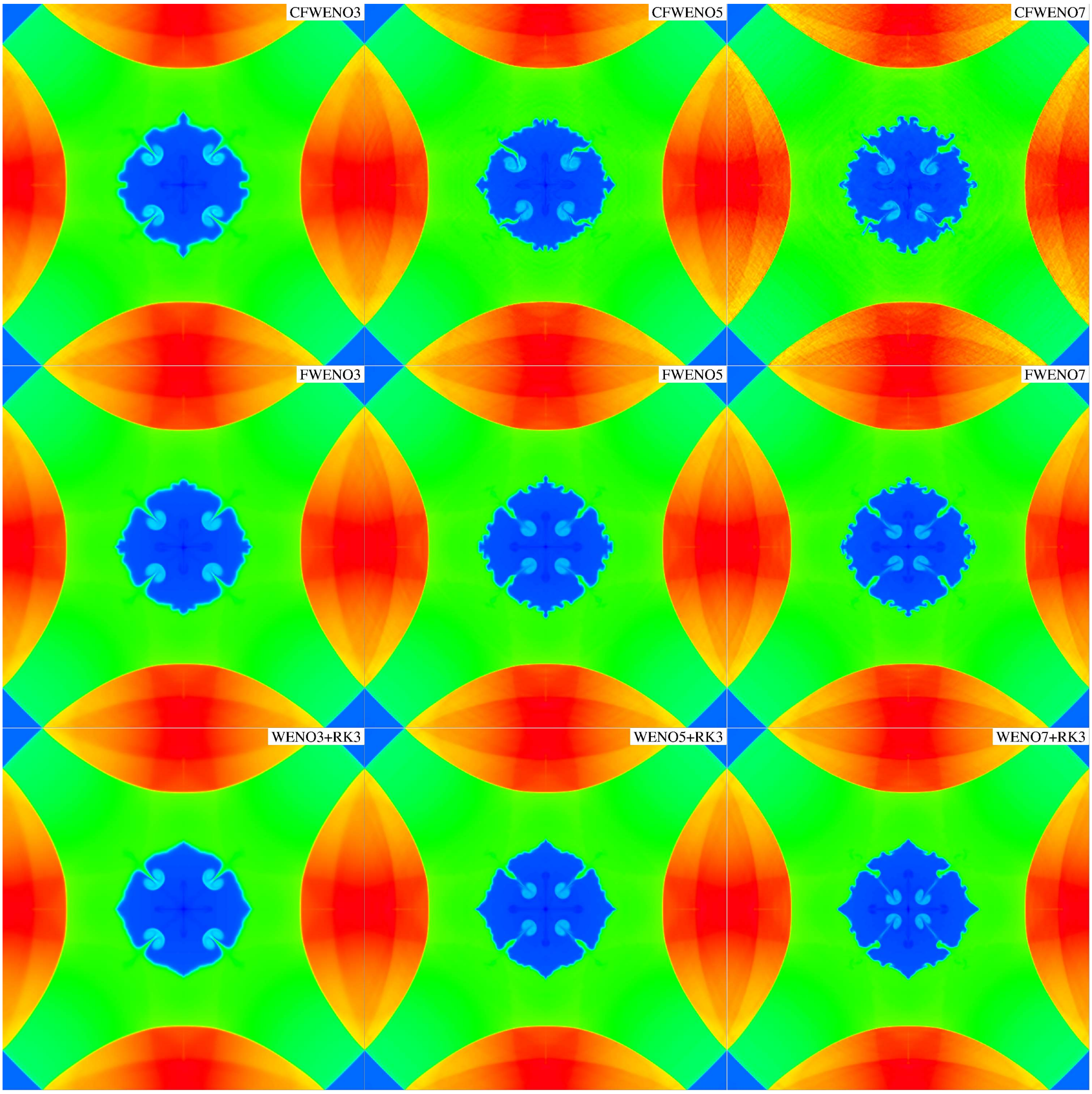

**Fig. 19 Two-dimensional Riemann problem: density contours of CFWENO, FWENO, WENO+RK3 (from top to below) and 3th,5th,7th order schemes (from left to right). Resolution 600×600. The simulation is run until $t$ = 0.4.**

**Example 4.4.3 Single-Material Triple point problem [64]**

A modified triple point problem with a single material rather than multiple materials is presented [64]. The computational domain is [0, 7]×[0, 3] and the initial condition is shown in Fig. 20. Outflow condition is applied to all boundaries. A uniform mesh with the resolution of 560×240 and 1120 × 480 is employed respectively for computations in Fig. 21 and Fig. 22.

Fig. 21 and Fig. 22 shows that at resolution 560 × 240, each third-order and fifth-order schemes perform similarly. For seventh-order schemes, CFWENO7 captures the most fine structures. Especially the slip line in the middle, only CFWENO5 and CFWENO7 roll up some small vortices. At the resolution of 1120 × 480, each scheme resolves significantly better than that of 560 × 240 grids. In the third-order scheme, CFWENO3 can capture some vortex structure on the interface, while FWENO3 and WENO3+RK3 basically cannot realize. The fifth-order and seventh-order have been significantly improved compared to the third-order. Among them, CFWENO has the highest resolution and captures more fine interface vortex structures. Especially for the middle slip line, CFWENO5 and CFWENO7 have captured obvious fine structures, while FWENO and WENO+RK3 have only seventh-order captured some of that and are not obvious.

In terms of efficiency, due to the special evolution method of two-dimensional CFWENO, the total number of grids refers to the sum of all node and half points, while FWENO and WENO+RK3 represent for only node points. The efficiency test in section 4.5 shows that the computational of CFWENO is about 1/3 of FWENO, about 1/10 to 1/13 of WENO+RK3; if the number of CFWENO grids is doubled, and keep the spatial step length in the evolution direction same as that of FWENO and WENO+RK3, the computational cost of CFWENO is about 2.6 to 3.0 times that of FWENO and about 0.6 to 0.8 times that of WENO+RK3. This reflects the extremely high efficiency of CFWENO.

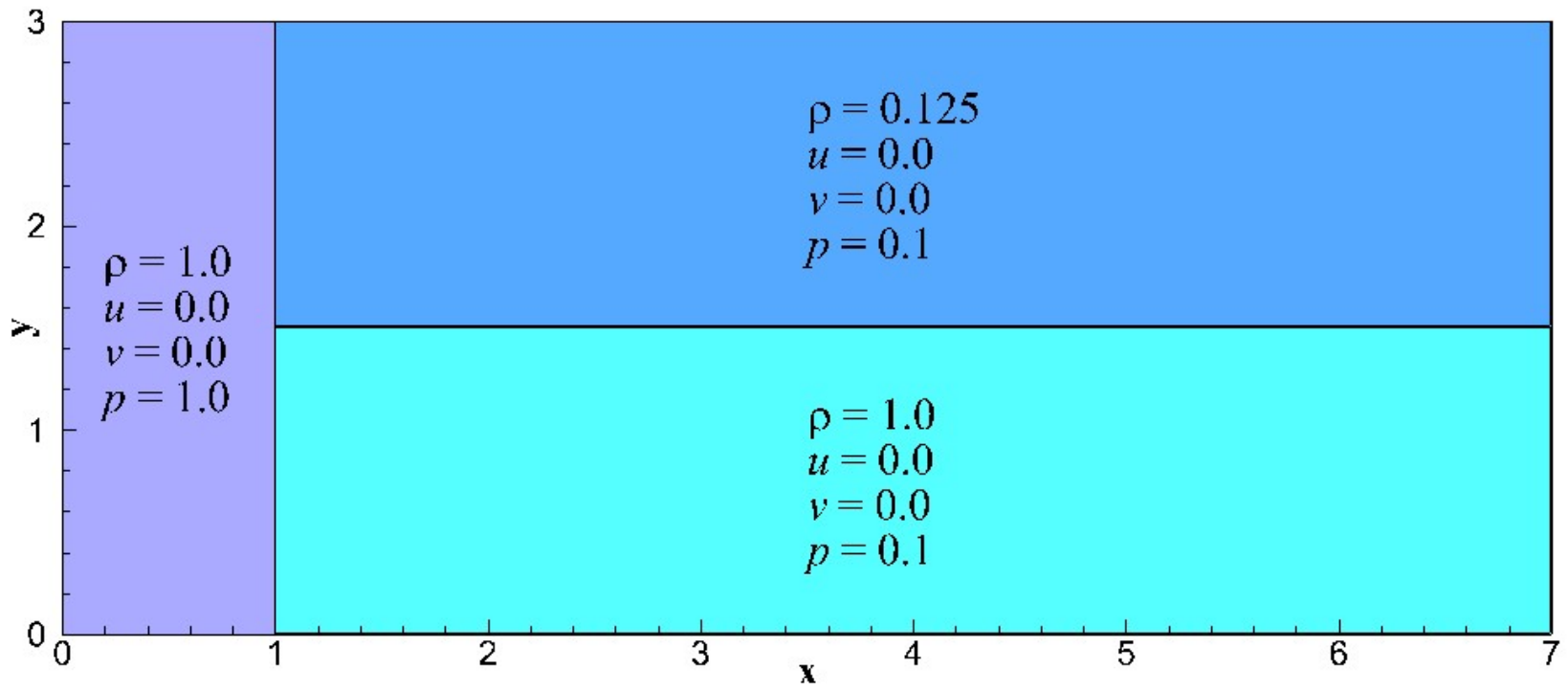

**Fig. 20 The sketch of the computational domain and the initial condition for the triple point problem**

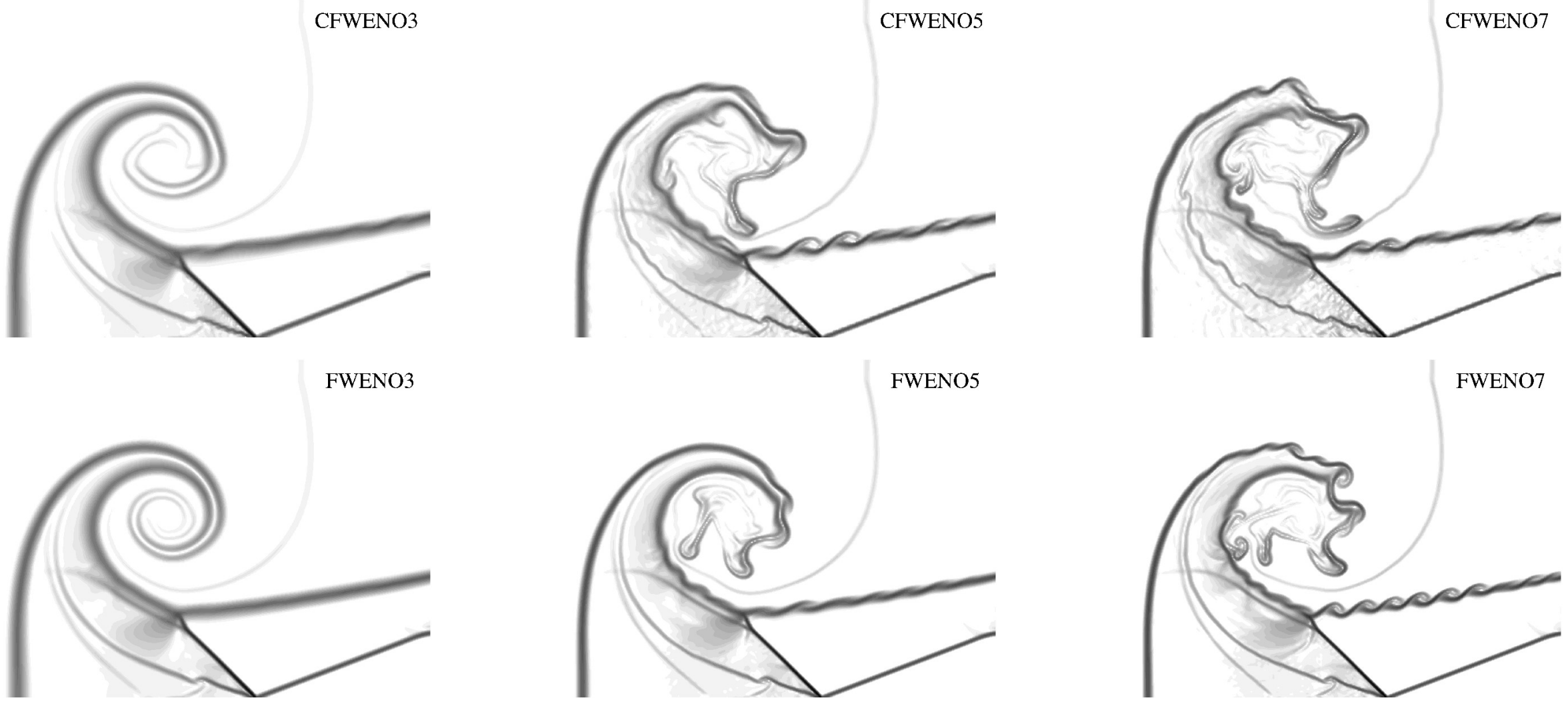

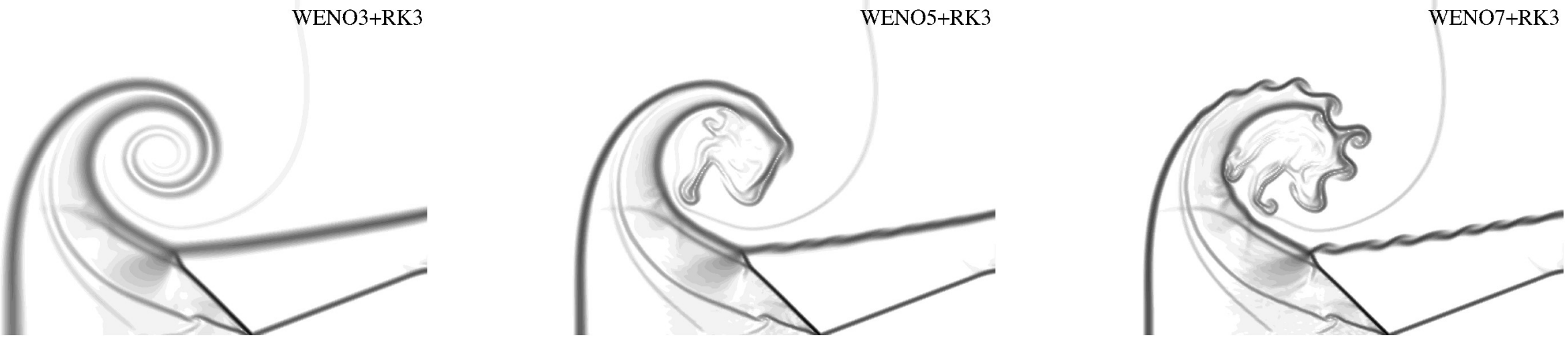

**Fig. 21 Single-Material Triple point problem resolution: normalized density gradient magnitude of CFWENO, FWENO, WENO+RK3 (from top to below) and 3th, 5th, 7th order schemes (from left to right). Resolution 560×240. The simulation is run until *t* = 5.0**

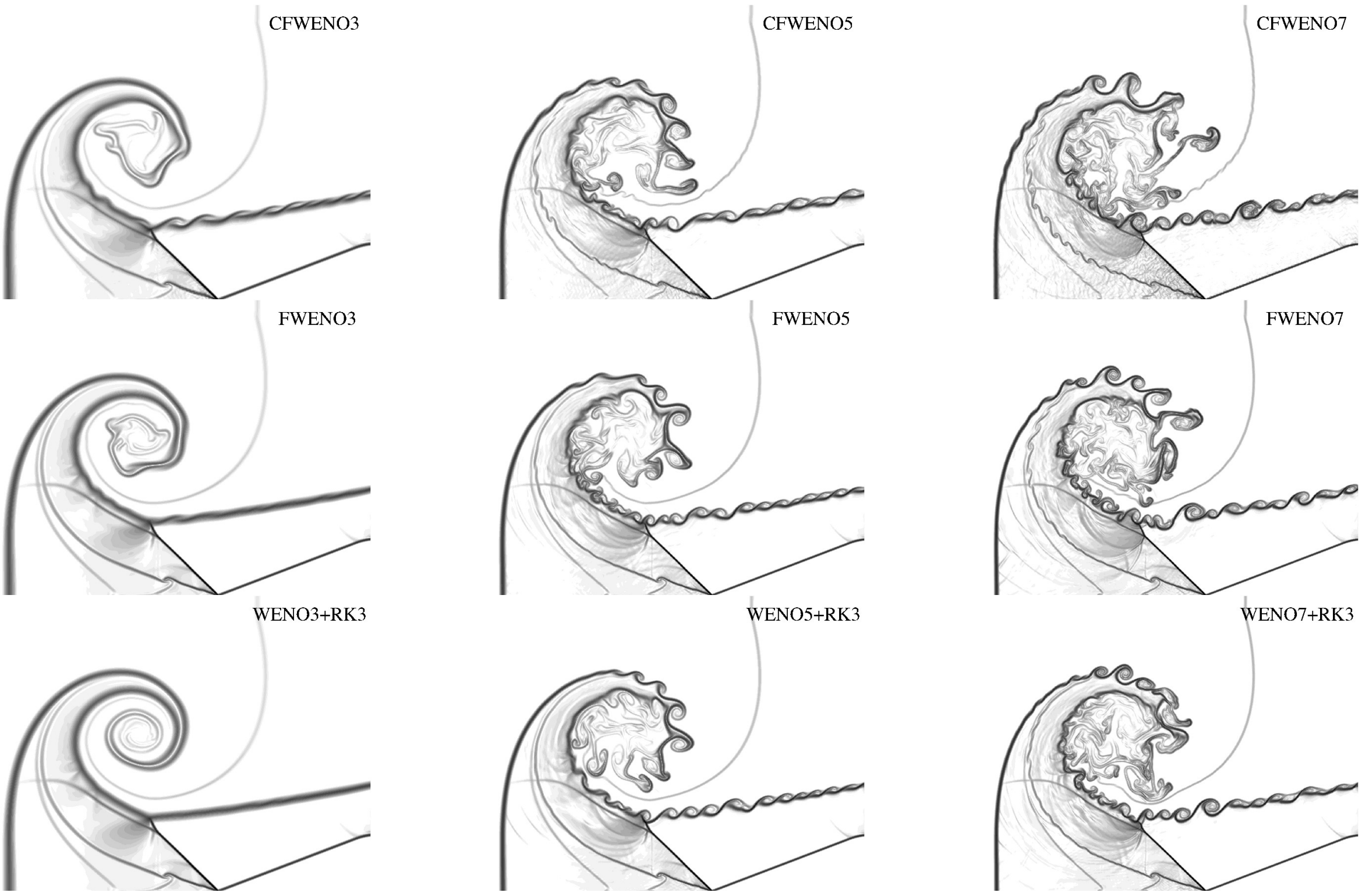

**Fig. 22 Single-Material Triple point problem: normalized density gradient magnitude of CFWENO, FWENO, WENO+RK3 (from top to below) and 3th, 5th, 7th order schemes (from left to right). Resolution 1120×480. The simulation is run until *t* = 5.0**

## 4.5 Efficiency test

In this section, the efficiency test is carried out. The test includes the scalar equation, one-dimensional and two-dimensional Euler equations are the examples in section 4.1 to 4.3, in which CFWENO and FWENO use CFL = 0.9, and WENO+RK3 uses the commonly used CFL = 0.6. In the test, the grid number of one-dimensional cases refers to the number of node values $u_i$, and the half point $u_{i+1/2}$ is not included. For the two-dimensional cases, the number of grids refers to the sum of all node and half points due to the newly designed evolution method for CFWENO in multi-dimensional problems. Table 9, Table 10, Table 11 show the computational cost of 3th, 5th and 7th order schemes respectively. It can be observed that:

(1) For the one-dimensional scalar equation, with the same grid resolution, the computational cost of CFWENO

is about 1.2 to 1.4 times that of FWENO, and about 0.25 times that of WENO+RK3.

(2) For the one-dimensional Euler equations, with the same grid resolution, the computational cost of CFWENO is about 1.2 to 1.4 times that of FWENO, and about 0.3~0.4 times that of WENO+RK3.

(3) For the two-dimensional Euler equations, with the same grid resolution, the computational cost of CFWENO is about 1/3 times that of FWENO, and about 1/10 to 1/13 times that of WENO+RK3, which is consistent with the theoretical estimate in Table 8. At the same time, the resolution of CFWNEO is similar to FWENO and WENO+RK3, even show higher resolution in many cases.

(4) For the two-dimensional Euler equation, if the number of grids in each direction of CFWENO is doubled, FWENO and WENO remain unchanged (that is, the spatial step length of CFWENO in the evolution direction is same as that of FWENO and WENO), the computational cost of CFWENO is about 2.6 to 3.0 times that of FWENO, and about 0.6 to 0.8 times that of WENO+RK3, which is also consistent with the inference in Table 5.

In summary, CFWENO is much more efficient than FWENO and WENO + RK.

**Table 9 Computational cost of 3th order schemes**

| Test example | Simulation time | Grid size | Wall-clock time (s) | | |
|---|---|---|---|---|---|
| | | | CFWENO3 | FWENO3 | WENO3+RK3 |
| **1D Scalar Equation** | | | | | |
| Burgers equation | 20.00 | 2000 | 4.98 | 4.15 | 19.50 |
| **1D Euler Equation** | | | | | |
| Sod | 0.20 | 2000 | 1.76 | 1.25 | 4.42 |
| Shu-Osher | 1.80 | 2000 | 3.26 | 2.61 | 9.05 |
| **2D Euler Equation** | | | | | |
| 2D Riemann 1 | 0.20 | 400×400 | 84.43 | 251.40 | 888.27 |
| | | 800×800 | 681.31 | —— | —— |
| Single-Material Triple point problem | 5.00 | 560×240 | 156.20 | 445.56 | 1485.63 |
| | | 1120×480 | 1184.11 | —— | —— |

**Table 10 Computational cost of 5th order schemes**

| Test example | Simulation time | Grid size | Wall-clock time (s) | | |
|---|---|---|---|---|---|
| | | | CFWENO5 | FWENO5 | WENO5+RK3 |
| **1D Scalar Equation** | | | | | |
| Burgers equation | 20.00 | 2000 | 8.66 | 6.40 | 34.27 |
| **1D Euler Equation** | | | | | |
| Sod | 0.20 | 2000 | 2.45 | 1.73 | 7.37 |
| Shu-Osher | 1.80 | 2000 | 4.64 | 3.71 | 14.51 |
| **2D Euler Equation** | | | | | |
| 2D Riemann 1 | 0.20 | 400×400 | 107.91 | 323.07 | 1288.45 |
| | | 800×800 | 968.21 | —— | —— |
| Single-Material Triple point problem | 5.00 | 560×240 | 202.21 | 618.59 | 2149.82 |
| | | 1120×480 | 1755.94 | —— | —— |

**Table 11 Computational cost of 7th order schemes**

| Test example | Simulation time | Grid size | Wall-clock time (s) | | |
|---|---|---|---|---|---|
| | | | CFWENO7 | FWENO7 | WENO7+RK3 |
| **1D Scalar Equation** | | | | | |
| Burgers equation | 20.00 | 2000 | 13.60 | 10.07 | 57.95 |
| **1D Euler Equation** | | | | | |
| Sod | 0.20 | 2000 | 3.56 | 2.51 | 11.63 |
| Shu-Osher | 1.80 | 2000 | 6.74 | 5.25 | 24.56 |
| **2D Euler Equation** | | | | | |
| 2D Riemann 1 | 0.20 | 400×400 | 141.11 | 423.79 | 1901.66 |
| | | 800×800 | 1164.56 | —— | —— |
| Single-Material Triple point problem | 5.00 | 560×240 | 280.89 | 828.84 | 3347.53 |
| | | 1120×480 | 2387.35 | —— | —— |

## 5 Conclusions

The solution formula method is a straightforward and clear approach for constructing one-step fully-discrete schemes, which is also naturally upwind and consistent high-order in both space and time. Under the framework of solution formula method, conservative variables at the half point (grid center) can be also obtained. We store these quantities and use them as interpolation information for the construction of high-order stencils, resulting in a new class of compact fully-discrete schemes. These schemes can be implemented with most stencil strategies. Using WENO reconstruction as an example, this paper proposes the CFWENO scheme. We also design a new entropy condition flux linearization reconstruction strategy for Euler equation, which enhances the robustness of CFWENO. In this paper, we conduct a detailed analysis:

(1) The accuracy analysis indicates that its nonlinear accuracy differs from that of FWENO. In CFWENO, half point values are used not only for flux reconstruction but also for initial value reconstruction. Therefore, each iteration of CFWENO improves the accuracy by only one order, while FWENO, used only for flux reconstruction, improves by two orders. Fortunately, experiments show that in practical computations, iterations are generally not necessary, and their impact on the results is not significant.

(2) Error analysis indicates that the error coefficient of CFWENO is significantly smaller than that of FWENO and WENO. Unlike FWENO, whose error coefficient monotonically decreases, and WENO, whose error coefficient remains constant, the error coefficient of CFWENO first increases and then decreases with the CFL number. Moreover, the error coefficient of CFWENO tends to be zero as the CFL number approaches either 0 or 1, while FWENO only tends to be zero when CFL is close to 1.

(3) Although CFWENO belongs to the Hermite interpolation scheme, it does not need to add any derivative equation because the HJ equation and the conservation law are used. Even there are only the conservative variables after simplification. However, the semi-discrete HWENO with RK method adopts the conservation law and its derivative equation, and the derivative equation also needs to be discretized. Also, there are cross derivatives in the multi-dimensional case that need to be handled for HWENO. The accuracy of stencils that CFWENO and HWENO can construct in the same dependency domain is also different. In $[x_{i-k}, x_{i+k}]$, CFWENO can construct $(4k+1)$-th order stencil at most, while HWENO can reach $(4k+2)$-th order. Nevertheless, in $[x_{i-(k+1/2)}, x_{i+(k+1/2)}]$, CFWENO can construct $(4k+3)$-th order stencil at most, while HWENO still only can reach $(4k+2)$-th order.

(4) The experimental results show that for the one-dimensional Euler equation, the calculation time of CFWENO is about 1.2 to 1.4 times that of FWENO, and about 0.3 to 0.4 times that of WENO+RK3. For the multi-dimensional case, we design a new evolution method. Although the new method increases the computing burden, the efficiency is still very high. For the two-dimensional Euler equations, due to the special evolution method of CFWENO, the half points are also counted in total grid resolution. Experiments show that, with the same grid resolution, the computational cost of CFWENO is about 1/3 that of FWENO and about 1/10 to 1/13 that of WENO+RK3. The conclusion remains consistent with theoretical analysis. In this situation, the resolution of CFWENO is still comparable to or even higher than that of FWENO and WENO+RK3, which indicates the high-efficiency property of CFWENO.

The high-resolution and high-efficiency properties of this compact fully-discrete framework give it great potential in unsteady computations such as turbulence. In the future, we will develop more robust schemes based on this framework that can handle various complex situations.

# Appendix A: Basic stencils for CFWENO

For $\overline{u}^{r}_{i+1/2,k}(k=0,1,\cdots,r-1)$, we have

$$r=2\quad\begin{cases}\overline{u}^{2}_{i+1/2,0}=u_j+(1-v)(u_j-u_{j-\frac{1}{2}})\\ \overline{u}^{2}_{i+1/2,1}=u_j+(1-v)(u_{j+\frac{1}{2}}-u_j)\\ \overline{u}^{2}_{i+1/2}\;=u_j+(1-v)(u_{j+\frac{1}{2}}-u_j)+(1-v)(-v)(u_{j-\frac{1}{2}}-2u_j+u_{j+\frac{1}{2}})\end{cases}\tag{A1}$$

$$r=3\quad\begin{cases}\overline{u}^{3}_{i+1/2,0}=u_j+(1-v)(u_j-u_{j-\frac{1}{2}})+\frac{1}{2}(1-v)^2(u_{j-1}-2u_{j-\frac{1}{2}}+u_j)\\ \overline{u}^{3}_{i+1/2,1}=u_j+(1-v)(u_{j+\frac{1}{2}}-u_j)+(1-v)(-v)(u_{j-\frac{1}{2}}-2u_j+u_{j+\frac{1}{2}})\\ \overline{u}^{3}_{i+1/2,2}=u_j+(1-v)(u_{j+\frac{1}{2}}-u_j)+\frac{1}{2}(1-v)(-v)(u_j-2u_{j+\frac{1}{2}}+u_{j+1})\\ \overline{u}^{3}_{i+1/2}\;=u_j+(1-v)(u_{j+\frac{1}{2}}-u_j)+\frac{1}{2}(1-v)(-v)(u_j-2u_{j+\frac{1}{2}}+u_{j+1})\\ \qquad+\frac{1}{4}(1-v)(-v)(1+v)(2u_{j-\frac{1}{2}}-5u_j+4u_{j+\frac{1}{2}}-u_{j+1})\\ \qquad+\frac{1}{12}(1-v)^2(-v)(1+v)(-u_{j-1}+6u_{j-\frac{1}{2}}-10u_j+6u_{j+\frac{1}{2}}-u_{j+1})\end{cases}\tag{A2}$$

$$r=4\quad\begin{cases}\overline{u}^{4}_{i+1/2,0}=u_j+(1-v)(u_j-u_{j-\frac{1}{2}})+\frac{1}{2}(1-v)^2(u_{j-1}-2u_{j-\frac{1}{2}}+u_j)\\ \qquad+\frac{1}{4}(2-v)(1-v)^2(-2u_{j-\frac{3}{2}}+5u_{j-1}-4u_{j-\frac{1}{2}}+u_j)\\ \overline{u}^{4}_{i+1/2,1}=u_j+(1-v)(u_{j+\frac{1}{2}}-u_j)+(1-v)(-v)(u_{j-\frac{1}{2}}-2u_j+u_{j+\frac{1}{2}})\\ \qquad+\frac{1}{4}(1-v)(-v)(1+v)(-u_{j-1}+4u_{j-\frac{1}{2}}-5u_j+2u_{j+\frac{1}{2}})\\ \overline{u}^{4}_{i+1/2,2}=u_j+(1-v)(u_{j+\frac{1}{2}}-u_j)+\frac{1}{2}(1-v)(-v)(u_j-2u_{j+\frac{1}{2}}+u_{j+1})\\ \qquad+\frac{1}{4}(1-v)(-v)(1+v)(2u_{j-\frac{1}{2}}-5u_j+4u_{j+\frac{1}{2}}-u_{j+1})\\ \overline{u}^{4}_{i+1/2,3}=u_j+(1-v)(u_{j+\frac{1}{2}}-u_j)+\frac{1}{2}(1-v)(-v)(u_j-2u_{j+\frac{1}{2}}+u_{j+1})\\ \qquad+\frac{1}{4}(1-v)(-v)(1+v)(u_j-4u_{j+\frac{1}{2}}+5u_{j+1}-2u_{j+\frac{3}{2}})\\ \overline{u}^{4}_{i+1/2}\;=u_j+(1-v)(u_{j+\frac{1}{2}}-u_j)+\frac{1}{2}(1-v)(-v)(u_j-2u_{j+\frac{1}{2}}+u_{j+1})\\ \qquad+\frac{1}{4}(1-v)(-v)(1+v)(2u_{j-\frac{1}{2}}-5u_j+4u_{j+\frac{1}{2}}-u_{j+1})\\ \qquad+\frac{1}{12}(1-v)^2(-v)(1+v)(-u_{j-1}+6u_{j-\frac{1}{2}}-10u_j+6u_{j+\frac{1}{2}}-u_{j+1})\\ \qquad+\frac{1}{36}(2-v)(1-v)^2(-v)(1+v)(3u_{j-\frac{3}{2}}-10u_{j-1}+18u_{j-\frac{1}{2}}-19u_j+9u_{j+\frac{1}{2}}-u_{j+1})\\ \qquad+\frac{1}{108}(2-v)^2(1-v)^2(-v)(1+v)\\ \qquad\quad(-3u_{j-\frac{3}{2}}+11u_{j-1}-27u_{j-\frac{1}{2}}+38u_j-27u_{j+\frac{1}{2}}+11u_{j+1}-3u_{j+\frac{3}{2}})\end{cases}\tag{A3}$$

For $u_{i+1/2,k}^{n+1,r}(k=0,1,\cdots,r-1)$, we have

$$r=2\quad\begin{cases} u_{i+1/2,0}^{n+1,2}=u_j+\dfrac{d}{dv}\big(-(1-v)(-v)\big)(u_j-u_{j-\frac{1}{2}}) \\ u_{i+1/2,1}^{n+1,2}=u_j+\dfrac{d}{dv}\big(-(1-v)(-v)\big)(u_{j+\frac{1}{2}}-u_j) \\ u_{i+1/2}^{n+1,2}=u_j+\dfrac{d}{dv}\big(-(1-v)(-v)\big)(u_{j+\frac{1}{2}}-u_j)+\dfrac{d}{dv}\Big(-\dfrac{1}{2}(1-v)(-v)^2\Big)(u_j-2u_{j+\frac{1}{2}}+u_{j+1}) \end{cases}\tag{A4}$$

$$r=3\quad\begin{cases} u_{i+1/2,0}^{n+1,3}=u_j+\dfrac{d}{dv}\big(-(1-v)(-v)\big)(u_j-u_{j-\frac{1}{2}})+\dfrac{d}{dv}\Big(-\dfrac{1}{2}(1-v)^2(-v)\Big)(u_{j-1}-2u_{j-\frac{1}{2}}+u_j) \\ u_{i+1/2,1}^{n+1,3}=u_j+\dfrac{d}{dv}\big(-(1-v)(-v)\big)(u_{j+\frac{1}{2}}-u_j)+\dfrac{d}{dv}\big(-(1-v)(-v)^2\big)(u_{j-\frac{1}{2}}-2u_j+u_{j+\frac{1}{2}}) \\ u_{i+1/2,2}^{n+1,3}=u_j+\dfrac{d}{dv}\big(-(1-v)(-v)\big)(u_{j+\frac{1}{2}}-u_j)+\dfrac{d}{dv}\Big(-\dfrac{1}{2}(1-v)(-v)^2\Big)(u_j-2u_{j+\frac{1}{2}}+u_{j+1}) \\ u_{i+1/2}^{n+1,3}=u_j+\dfrac{d}{dv}\big(-(1-v)(-v)\big)(u_{j+\frac{1}{2}}-u_j)+\dfrac{d}{dv}\Big(-\dfrac{1}{2}(1-v)(-v)^2\Big)(u_j-2u_{j+\frac{1}{2}}+u_{j+1}) \\ \qquad+\dfrac{d}{dv}\Big(-\dfrac{1}{4}(1-v)(-v)^2(1+v)\Big)(2u_{j-\frac{1}{2}}-5u_j+4u_{j+\frac{1}{2}}-u_{j+1}) \\ \qquad+\dfrac{d}{dv}\Big(-\dfrac{1}{12}(1-v)^2(-v)^2(1+v)\Big)(-u_{j-1}+6u_{j-\frac{1}{2}}-10u_j+6u_{j+\frac{1}{2}}-u_{j+1}) \end{cases}\tag{A5}$$

$$r=4\quad\begin{cases} u_{i+1/2,0}^{n+1,4}=u_j+\dfrac{d}{dv}\big(-(1-v)(-v)\big)(u_j-u_{j-\frac{1}{2}})+\dfrac{d}{dv}\Big(-\dfrac{1}{2}(1-v)^2(-v)\Big)(u_{j-1}-2u_{j-\frac{1}{2}}+u_j) \\ \qquad+\dfrac{d}{dv}\Big(-\dfrac{1}{4}(2-v)(1-v)^2(-v)\Big)(-2u_{j-\frac{3}{2}}+5u_{j-1}-4u_{j-\frac{1}{2}}+u_j) \\ u_{i+1/2,1}^{n+1,4}=u_j+\dfrac{d}{dv}\big(-(1-v)(-v)\big)(u_{j+\frac{1}{2}}-u_j)+\dfrac{d}{dv}\big(-(1-v)(-v)^2\big)(u_{j-\frac{1}{2}}-2u_j+u_{j+\frac{1}{2}}) \\ \qquad+\dfrac{d}{dv}\Big(-\dfrac{1}{4}(1-v)(-v)^2(1+v)\Big)(-u_{j-1}+4u_{j-\frac{1}{2}}-5u_j+2u_{j+\frac{1}{2}}) \\ u_{i+1/2,2}^{n+1,4}=u_j+\dfrac{d}{dv}\big(-(1-v)(-v)\big)(u_{j+\frac{1}{2}}-u_j)+\dfrac{d}{dv}\Big(-\dfrac{1}{2}(1-v)(-v)^2\Big)(u_j-2u_{j+\frac{1}{2}}+u_{j+1}) \\ \qquad+\dfrac{d}{dv}\Big(-\dfrac{1}{4}(1-v)(-v)^2(1+v)\Big)(2u_{j-\frac{1}{2}}-5u_j+4u_{j+\frac{1}{2}}-u_{j+1}) \\ u_{i+1/2,3}^{n+1,4}=u_j+\dfrac{d}{dv}\big(-(1-v)(-v)\big)(u_{j+\frac{1}{2}}-u_j)+\dfrac{d}{dv}\Big(-\dfrac{1}{2}(1-v)(-v)^2\Big)(u_j-2u_{j+\frac{1}{2}}+u_{j+1}) \\ \qquad+\dfrac{d}{dv}\Big(-\dfrac{1}{4}(1-v)(-v)^2(1+v)\Big)(u_j-4u_{j+\frac{1}{2}}+5u_{j+1}-2u_{j+\frac{3}{2}}) \\ u_{i+1/2}^{n+1,4}=u_j+\dfrac{d}{dv}\big(-(1-v)(-v)\big)(u_{j+\frac{1}{2}}-u_j)+\dfrac{d}{dv}\Big(-\dfrac{1}{2}(1-v)(-v)^2\Big)(u_j-2u_{j+\frac{1}{2}}+u_{j+1}) \\ \qquad+\dfrac{d}{dv}\Big(-\dfrac{1}{4}(1-v)(-v)^2(1+v)\Big)(2u_{j-\frac{1}{2}}-5u_j+4u_{j+\frac{1}{2}}-u_{j+1}) \\ \qquad+\dfrac{d}{dv}\Big(-\dfrac{1}{12}(1-v)^2(-v)^2(1+v)\Big)(-u_{j-1}+6u_{j-\frac{1}{2}}-10u_j+6u_{j+\frac{1}{2}}-u_{j+1}) \\ \qquad+\dfrac{d}{dv}\Big(-\dfrac{1}{36}(2-v)(1-v)^2(-v)^2(1+v)\Big)(3u_{j-\frac{3}{2}}-10u_{j-1}+18u_{j-\frac{1}{2}}-19u_j+9u_{j+\frac{1}{2}}-u_{j+1}) \\ \qquad+\dfrac{d}{dv}\Big(-\dfrac{1}{108}(2-v)^2(1-v)^2(-v)^2(1+v)\Big) \\ \qquad(-3u_{j-\frac{3}{2}}+11u_{j-1}-27u_{j-\frac{1}{2}}+38u_j-27u_{j+\frac{1}{2}}+11u_{j+1}-3u_{j+\frac{3}{2}}) \end{cases}\tag{A6}$$